\documentclass[preprint,12pt]{article}
\usepackage{amssymb}
\usepackage{mathrsfs}
\usepackage{amsfonts}
\usepackage{graphicx}
\usepackage{colortbl,dcolumn}
\usepackage{amsmath}
\usepackage{psfrag}
\usepackage{booktabs}
\usepackage{slashbox}
\allowdisplaybreaks
\numberwithin{equation}{section}
\usepackage[top=1in, bottom=1in, left=0.9in, right=0.9in]{geometry}

\newtheorem{thm}{Theorem}[section]

\newtheorem{lem}[thm]{Lemma}

\newtheorem{ass}[thm]{Assumption}
\begin{document}
%\title{A higher-order explicit strongly convergent numerical method for
\title{The tamed  Milstein method for commutative
stochastic differential equations with non-globally Lipschitz continuous coefficients \footnotemark[1]}
       \author{
        Xiaojie Wang \quad Siqing Gan \footnotemark[2] \\ %\quad  Siqing Gan \footnotemark[3]
       {\small School of Mathematics and Statistics,
       Central South University,}\\
      {\small Changsha 410083, Hunan,  PR China } }

       \date{}
       \maketitle

       \footnotetext{\footnotemark[1] This work was supported by NSF of China (No.11171352) and Hunan Provincial Innovation Foundation For Postgraduate (NO.CX2010B118). The first author would like to thank Professor Peter Kloeden and Professor Arnulf Jentzen for their kindness during his stay in Frankfurt from September 2010 to August 2011.}
        \footnotetext{\footnotemark[2]Corresponding author: x.j.wang7@gmail.com}
       %\footnotetext{\footnotemark[3]Corresponding author: siqinggan@yahoo.com.cn}
       \begin{abstract}
          {\rm\small For stochastic differential equations (SDEs) with a superlinearly growing and globally one-sided Lipschitz continuous drift coefficient, the classical explicit Euler scheme fails to converge strongly to the exact solution. Recently, an explicit strongly convergent numerical scheme, called the tamed Euler method, is proposed in [Hutzenthaler, Jentzen, $\&$ Kloeden, Ann. Appl. Probab., 22 (2012), pp. 1611-1641.] for such SDEs. Motivated by their work, we here introduce a tamed version of the Milstein scheme for SDEs with commutative noise. The proposed method is also explicit and easily implementable, but achieves higher strong convergence order than the tamed Euler method does. In recovering  the strong convergence order one of the new method, new difficulties arise and kind of a bootstrap argument is developed to overcome them. Finally, an illustrative example confirms the computational efficiency of the tamed Milstein method compared to the tamed Euler method. }\\

\textbf{AMS subject classification: } {\rm\small 65C20, 60H35, 65L20.}\\
%\textbf{PACS: 02.60.Lj}

\textbf{Key Words: }{\rm\small} tamed Milstein method, superlinearly growing coefficient, one-sided Lipschitz condition,
commutative noise, strong convergence
\end{abstract}

%% main text
\section{Introduction}
\label{introduction}

We consider numerical integration of stochastic differential equations (SDEs) in the It\^{o}'s sense
\begin{equation}\label{sdes}
dX_t=\mu(X_t)dt + \sigma(X_t)dW_t, \quad X_0= \xi, \quad t \in[0,T].
\end{equation}
Here $\mu:\mathbb{R}^d\longrightarrow\mathbb{R}^d,
\sigma= (\sigma_1, \sigma_2,...,\sigma_m):\mathbb{R}^d\longrightarrow\mathbb{R}^{d\times m}$. We assume that
$W_t$ is an $m$-dimensional Wiener process defined on the complete
probability space $(\Omega,\mathcal {F},\mathbb{P})$ with an increasing filtration $\{\mathcal {F}_t\}_{t\geq 0}$
satisfying the usual conditions. And the initial data $\xi$ is independent of the Wiener process. (\ref{sdes}) can be interpreted mathematically as a stochastic
integral equation
\begin{equation}\label{itsde}
X_t=X_0 + \int_0^t \mu(X_s)ds + \sum_{i=1}^{m}\int_0^t \sigma_i(X_s)dW_s^i, \quad t
\in[0,T], \mathbb{P}-a.s.,
\end{equation}
where $\sigma_i(x) = (\sigma_{1,i}(x),..., \sigma_{d,i}(x))^T$ for $x \in \mathbb{R}^d, i \in \{1,2,...,m\}$ and the second integral is the It\^{o} integral.

This article is concerned with the strong approximation problem (see, e.g., Section 9.3 in Kloeden and Platen \cite{KP92}) of the SDEs (\ref{itsde}). More precisely, on a uniform mesh with stepsize $h = \frac{T}{N}$ defined by $ \mathcal {T}^N : {0 = t_0 < t_1 < t_2 < \cdots < t_N=T}, N \in \mathbb{N}$, we want to compute a numerical approximation $Y_n: \Omega \rightarrow \mathbb{R}^d, n\in \{0,1,...,N\}$ with $Y_0=\xi$ such that
\begin{equation} \label{Err}
\left(\mathbb{E} \|X_T-Y_N\|^2\right)^{\frac{1}{2}} < \varepsilon
\end{equation}
for a given precision $\varepsilon >0$ with the least possible computational effort.
%(number of arithmetical operations, i.e., addition, subtraction, multiplication and
%division, and independent standard normal random variables needed to compute $Y_N$).
%Now let us comment on the importance of solving the strong approximation problem
%(\ref{Err}). A central motivation for studying strong
%approximations in the sense of (\ref{Err})
%is Giles' seminal paper \cite{Giles}. There he introduces a very efficient
%accelerated Monte Carlo method, as opposed to the classical Monte Carlo method,
%for approximating moments or other expectations of functionals of the SDE solution
%via numerical schemes that converge strongly. In view of this method, strong approximation
%of the exact solution of the SDEs (\ref{itsde}) in the sense of (\ref{Err}) yield very
%efficient approximations of expectations of functionals of the SDE solution and
%this is a central reason for developing strongly convergent numerical methods.
The strong convergence problem  becomes very important because efficient Multi-Level Monte Carlo (MLMC) simulations rely on the strong convergence properties \cite{Giles}.

The simplest and most obvious idea to solve the strong approximation problem
(\ref{Err}) is to apply the explicit Euler scheme \cite{Maruyama55}
\begin{equation}  \label{Euler}
Y_{n+1}= Y_{n} + h\mu(Y_{n}) + \sigma (Y_{n})\Delta W_n, \quad Y_0 = \xi, \:n=0,1,...,N-1,
\end{equation}
where $\Delta W_n = W_{t_{n+1}}-W_{t_{n}}$. The method is strongly convergent with order one half if the coefficients $\mu,\sigma$ satisfy the global Lipschitz condition (see, for instance, \cite{KP92}). Unfortunately, it
has recently been shown in \cite{HJ09a} that the explicit Euler scheme fails to provide strong convergent solution to the SDEs with super-linearly growing drift coefficient. It is well-known that the backward Euler method can promise strong convergence in this situation, see e.g.,\cite{HMS02}. But the backward Euler method is an implicit method, which requires additional computational effort to solve an implicit system.
Recently in \cite{HJK10}, the authors proposed an explicit method, called tamed Euler method, for (\ref{itsde})
\begin{equation}  \label{tamedEuler}
\begin{split}
Y_{n+1}= &Y_{n} + h \tilde{\mu}(Y_n) + \sigma (Y_{n})\Delta
W_n, \quad \mbox{with} \quad \tilde{\mu}(Y_n) =\frac{\mu(Y_{n})}{1+h \|\mu(Y_{n})\|}.
\end{split}
\end{equation}
Here $\tilde{\mu}$ is a modification of $\mu$. This tamed Euler scheme is proved to converge strongly with the standard convergence order 0.5 to the exact solution
of (\ref{itsde}) if the drift coefficient function is globally one-sided Lipschitz continuous and has an at most polynomially growing derivative.

On the one hand, the explicit Milstein scheme is another  numerical
scheme for SDEs that achieves a strong order of
convergence higher than that of the explicit Euler scheme \eqref{Euler} \cite{JR10,KP92,GNM74}. In fact the explicit Milstein scheme has strong convergence order of one if the coefficient functions in the stochastic Taylor expansions satisfy both the global Lipschitz condition and the linear growth condition(see \cite{KP92}).  The explicit Milstein method \cite{KP92,GNM74} applied to \eqref{sdes} reads
\begin{equation}  \label{Milstein}
\begin{split}
Y_{n+1}= &Y_{n} + h\mu(Y_{n}) + \sigma (Y_{n})\Delta
W_n + \sum^m_{j_1,j_2=1}L^{j_1}\sigma_{j_2}(Y_{n})I_{j_1,j_2}^{t_n,t_{n+1}},
\end{split}
\end{equation}
where
\begin{equation}\label{Lj}
L^{j_1} = \sum^d_{k=1} \sigma_{k,j_1}\frac{\partial}{\partial x^k}, \quad
I_{j_1,j_2}^{t_n,t_{n+1}} = \int^{t_{n+1}}_{t_n} \int^{s_2}_{t_n} d W_{s_1}^{j_1}
d W_{s_2}^{j_2}.
\end{equation}
Since the explicit Milstein scheme and the explicit Euler scheme coincide when applied to the SDEs with additive noise, we can deduce from the results in \cite{HJ09a} that the explicit Milstein scheme generally does not converge in the mean-square sense to the exact solution solution of the SDEs with super-linearly growing drift coefficient.
Accordingly, we follow the idea from \cite{HJK10} and replace $\mu(Y_{n})$ in (\ref{Milstein}) with $\tilde{\mu}(Y_{n})$ to derive a tamed Milstein method
\begin{equation}  \label{TamedMilstein}
\begin{split}
Y_{n+1}= &Y_{n} + h \tilde{\mu}(Y_n) + \sigma (Y_{n})\Delta
W_n + \sum^m_{j_1,j_2=1}L^{j_1}\sigma_{j_2}(Y_{n})I_{j_1,j_2}^{t_n,t_{n+1}},
\end{split}
\end{equation}
which we expect to be strongly convergent with order one in the non-globally Lipschitz case.

On the other hand, although Milstein-type schemes may achieve a strong convergence order higher than that of Euler-type schemes, additional computational effort is required to approximate the iterated It\^{o} integrals $I_{j_1,j_2}^{t_n,t_{n+1}}$ for every time step \cite{KPW92}. This will enable the  Milstein-type schemes to lose their advantage over the Euler-type schemes in computational efficiency. In this article we restrict our attention to SDEs with commutative noise, in which case the Milstein scheme can be easily implemented without simulating the iterated It\^{o} integrals. In this situation, Milstein-type method is much more computationally efficient than Euler-type method. More precisely, let the diffusion matrix $\sigma$ fulfill the so-called commutativity condition:
\begin{equation} \label{cc}
L^{j_1}\sigma_{k,j_2} = L^{j_2}\sigma_{k,j_1}, \quad j_1,j_2 = 1,...,m, k = 1,...,d.
\end{equation}
In many applications the considered SDE systems possess commutative noise (see \cite{KP92}).

%For example, SDEs with additive noise, scalar noise, diagonal noise in case of
%$d=m$ with $\sigma_{i,j}(x) \equiv 0$ and $ \partial \sigma_{i,i}(x)/\partial
%x^j \equiv 0$ for $i \neq j$, and linear noise with $\sigma_{i,j}(x) =
%\sigma_{i,j}x^i$  for $x=(x_1,...,x_d)^T\in \mathbb{R}^d, \sigma_{i,j} \in
%\mathbb{R}$ all satisfy the commutativity condition.

Thanks to the property $I_{j_1,j_2}^{t_n,t_{n+1}}+I_{j_2,j_1}^{t_n,t_{n+1}} =\Delta W_n^{j_1}\Delta W_n^{j_2}, j_1 \neq j_2$, in this case the tamed Milstein method (\ref{TamedMilstein}) takes a simple form as
\begin{equation}  \label{MilsteinC}
\begin{split}
Y_{n+1}= &Y_{n} + h \tilde{\mu}(Y_n) + \sigma (Y_{n})\Delta
W_n + \frac{1}{2}\sum^m_{j_1,j_2=1}L^{j_1}\sigma_{j_2}(Y_{n})\left(\Delta W^{j_1}_n \Delta W^{j_2}_n-\delta_{j_1,j_2}h\right),
\end{split}
\end{equation}
where $\delta_{j_1,j_2} = 1$ for $j_1=j_2$ and
$\delta_{j_1,j_2} = 0$ for $j_1\neq j_2$, $\tilde{\mu}$ is the
 modification of $\mu$ as defined in \eqref{tamedEuler}.
%$$
%\delta_{j_1,j_2} = \left\{\begin{array}{l}
%1,    \ \ \ j_1=j_2,\\
%0, \ \ \   j_1\neq j_2.
%\end{array}\right.
%$$

The main result of this article shows that the tamed Milstein scheme (\ref{MilsteinC})
converges strongly with the standard convergence order one to the exact solution
of SDEs with commutative noise if the drift coefficient $\mu$ is globally one-sided Lipschitz continuous and has at most polynomially growing first and second derivatives. The diffusion coefficient  $\sigma$ and the coefficient function $L^{j_1}\sigma_{j_2}, j_1,j_2 \in \{1,...,m\}$ are assumed to be globally Lipschitz continuous. It is worthwhile to mention that a similar approach as used in \cite{HJK10} is evoked to obtain uniform boundedness of $p$-th moments of numerical solutions produced by the tamed Milstein method.  We also introduce similar stochastic processes $D_n$ that dominate the tamed Milstein approximation on appropriate subevents $\Omega_n$ (see Section 2 for more details). With bounded $p$-th moments at hand, our main effort is to show for the time continuous tamed Milstein method there exists a family of real numbers $C_{p,T} \in [1,\infty)$ for $p\in [1,\infty)$ such that
\begin{equation} \label{msc}
\Big(\mathbb{E} \Big[\sup_{t\in[0,T]}\Big\|X_t -\bar{Y}_t\Big\|^p\Big]\Big)^{1/p} \leq C_{p,T}\cdot h, \quad h \in (0, 1].
\end{equation}
The key difficulty is that relative to previous analysis \cite{HJK10} a sharper estimate of the term $J$ (see \eqref{eq:J}) must be obtained to get the strong convergence order one. To overcome this difficulty, a certain kind of bootstrap argument is exploited (see the estimate of $J$ for more details). To the best of our knowledge, this is the very first paper to successfully recover the strong convergence order one for the Milstein-type method under non-globally Lipschitz condition.

% Compared to the proof of tamed Euler method in \cite{HJK10}, new difficulties arise in
%recovering the strong convergence order one of the method \eqref{MilsteinC} in the term
%of \eqref{msc}. In this work new techniques are developed to overcome the difficulties.

%We also mention that in many applications the considered SDE system possess
%the commutative noise. For instance,
%\begin{equation} \label{Omega}
%\Omega_n := \{\omega \in\Omega| \sup_{0\leq k\leq n-1}D_k(\omega)\leq \frac{1}{N^{2c}},}
%|I_{j_1,j_2}^{t_k,t_{k+1}}|\leq 1 \}
%\sup_{0\leq k\leq n-1} \|\Delta W_k\|\leq 1, \sup_{0\leq k\leq n-1} \sum_{j_1,j_2=1}^{m
%\end{equation}
%\begin{equation} \label{Dn}
%\begin{split}
%D_n &:=  (\lambda + \|\xi\|)\sqrt{\sup_{0\leq u\leq n}\prod_{k=u}^{n-1}
%\left[1+\lambda h + \lambda \|\Delta W_k\|^2+2\alpha_k + (3K^2+K^2h^{-1})|
%\sum_{j_1,j_2=1}^{m}I_{j_1,j_2}^{t_k,t_{k+1}}|^2\right] }
%\end{split}
%\end{equation}
The rest of this paper is arranged as follows. In the next section, uniform boundedness of $p$-th moments are obtained. And then the strong convergence order of the tamed Milstein method is established in Section 3. Finally, an illustrative example confirms the strong convergence order of one and the computational efficiency of this scheme compared to the tamed Euler scheme.

%\section{It\^{o} and Stratonovich stochastic differential equations}
\section{Uniform boundedness of $p$-th moments}

Throughout this article, $N \in \mathbb{N}$ is the step number of the uniform mesh defined in the previous section. Moreover, we use the notation $\|x\|:=(|x_1|^2 + ...+|x_k|^2)^{\frac{1}{2}}$, $\langle x, y \rangle := x_1y_1+...+x_ky_k$ for all $x = (x_1, x_2,...,x_k), y = (y_1, y_2,...,y_k) \in \mathbb{R}^k, k \in \mathbb{N}$, and $\|A\|:= \sup_{x\in\mathbb{R}^l,\|x\|\leq1}\|Ax\|$ for all $A\in\mathbb{R}^{k\times l}, k,l \in \mathbb{N}$. Furthermore, we make the following assumptions.

\begin{ass} \label{LCMC}
Let $\mu(x)$ and $\sigma_i(x), i=1,...,m$ be continuously differentiable and there exist positive constants $K \geq 1$ and $c \geq 1$, such that $\forall
x,y \in \mathbb{R}^d$
\begin{eqnarray}
\langle x-y, \mu(x)-\mu(y)\rangle &\leq& K\|x-y\|^2, \label{OLC1}\\
\|\sigma(x) - \sigma(y)\| &\leq& K\|x-y\|,\label{OLC3}\\
\|L^{j_1}\sigma_{j_2}(x)-L^{j_1}\sigma_{j_2}(y)\|&\leq& K\|x-y\|, \quad j_1,j_2 \in \{1,2,...,m\}, \label{OLC4}\\
\|\mu'(x)\| &\leq& K(1+\|x\|^c).\label{OLC5}
\end{eqnarray}
\end{ass}
Note that the globally one-sided Lipschitz  condition (\ref{OLC1}) on the drift $\mu$  and  the  globally Lipschitz  condition (\ref{OLC3}) on the diffusion $\sigma$  have been widely used in the literatures, e.g., \cite{HMS02,HK07,YH96,HJ09a,HJ09b,HJK10}.

To prove uniform boundedness of $p$-th moments of the numerical solution, we follow the ideas in \cite{HJK10} to introduce the appropriate subevents $\Omega_n$ and dominating stochastic processes $D_n$
\begin{equation} \label{Omega}
\Omega_n := \Big\{\omega \in\Omega| \sup_{0\leq k\leq n-1}D_k(\omega)\leq \frac{1}{N^{2c}}, \sup_{0\leq k\leq n-1} \|\Delta W_k\|\leq 1 \Big\},
\end{equation}
\begin{equation} \label{Dn}
\begin{split}
D_n &:=  (\lambda + \|\xi\|)\exp\Big(\lambda+ \sup_{0\leq u\leq n}\sum_{k=u}^{n-1}\left[\lambda \|\Delta W_k\|^2+\alpha_k\right]\Big),
\end{split}
\end{equation}
where
$$
\lambda = \left(1+ 4TK + 2T\|\mu(0)\| + 2K + 2\|\sigma(0)\|+ m^2(T+1)\Big(K+\max_{1\leq j_1,j_2\leq m}\|L^{j_1}\sigma_{j_2}(0)\|\Big)\right)^2
$$
and
\begin{equation} \label{alpha}
\begin{split}
\alpha_n =  1_{\{\|Y_n\| \geq 1\}}\Big\langle\frac{Y_n}{\|Y_n\|}, \frac{\sigma(Y_n)}{\|Y_n\|}\Delta W_n\Big\rangle
+ 1_{\{\|Y_n\| \geq 1\}}\Big\langle\frac{Y_n}{\|Y_n\|}, \sum_{j_1,j_2=1}^m\frac{L^{j_1}\sigma_{j_2}(Y_n)}{2\|Y_n\|}\left(\Delta W^{j_1}_n \Delta W^{j_2}_n-\delta_{j_1,j_2}h\right)\Big\rangle.
\end{split}
\end{equation}
The following lemmas are needed in order to prove uniform boundedness of $p$-th moments.
\begin{lem} \label{DomLem}
Let $Y_n, D_n$ and $\Omega_n$ be given by (\ref{MilsteinC}),(\ref{Dn}) and (\ref{Omega}), respectively. Then
\begin{equation} \label{d1}
1_{\Omega_n}\|Y_n\| \leq D_n, \quad \mbox{for all} \:\: n =0,1,...,N.
\end{equation}
\end{lem}
{\it Proof.} First of all, note that $\|\Delta W_n\|\leq 1$  on $\Omega_{n+1}$ for all $0 \leq n\leq N-1$ and $N\in \mathbb{N}$. The globally Lipschitz continuity of $\sigma$ and $L^{j_1}\sigma_{j_2}, j_1,j_2 \in \{1,2,...,m\}$, and the polynomial growth bound on $\mu'$ imply that, on $\Omega_{n+1} \cap \{\omega \in \Omega|\|Y_{n}(\omega)\|\leq 1 \}$ for all $0 \leq n\leq N-1$,
\begin{equation}\label{yn1}
\begin{split}
&\|Y_{n+1}\| \\
\leq& \|Y_{n}\| + h \|\mu(Y_n)\| + \|\sigma(Y_n)\|\|\Delta W_n\| + \frac{1}{2} \sum_{j_1,j_2=1}^{m}\|L^{j_1}\sigma_{j_2}(Y_{n})\||\Delta W^{j_1}_n \Delta W^{j_2}_n-\delta_{j_1,j_2}h|
\\\leq& 1+ h\|\mu(Y_n)-\mu(0)\|+h\|\mu(0)\| + \|\sigma(Y_n)-\sigma(0)\|\|\Delta W_n\|+ \|\sigma(0)\|\|\Delta W_n\|
\\ &+ \frac{1}{2}\sum_{j_1,j_2=1}^{m}\left(\|L^{j_1}\sigma_{j_2}(Y_{n})- L^{j_1}\sigma_{j_2}(0)\|+\|L^{j_1}\sigma_{j_2}(0)\|\right)|\Delta W^{j_1}_n \Delta W^{j_2}_n-\delta_{j_1,j_2}h|
\\ \leq &1 + hK(1+\|Y_n\|^c)\|Y_n\|+ h \|\mu(0)\| + K\|Y_n\|\|\Delta W_n\|+\|\sigma(0)\|\|\Delta W_n\|
\\ &+ \frac{1}{2}\sum_{j_1,j_2=1}^{m}\left(K\|Y_{n}\|+\|L^{j_1}\sigma_{j_2}(0)\|\right)|\Delta W^{j_1}_n \Delta W^{j_2}_n-\delta_{j_1,j_2}h|
\\ \leq &1+ 2TK+T\|\mu(0)\|+K+\|\sigma(0)\| + \frac{m}{2}(m+T)\Big(K+\max_{1\leq j_1,j_2\leq m}\|L^{j_1}\sigma_{j_2}(0)\|\Big)
\\ \leq& \lambda.
\end{split}
\end{equation}
Moreover, the Cauchy-Schwarz inequality and the inequality $a\cdot b \leq \frac{a^2}{2}+\frac{b^2}{2}$ for all $a,b \in \mathbb{R}$ give that
\begin{equation}\label{es0}
\begin{split}
\|Y_{n+1}\|^2 =& \left\|Y_{n} + h \tilde{\mu}(Y_n) + \sigma (Y_{n})\Delta W_n +M_n\right\|^2
\\=& \|Y_n\|^2 + h^2\| \tilde{\mu}(Y_n)\|^2 + \|\sigma(Y_n)\Delta W_n\|^2 +\|M_n\|^2
\\&  +
2h \langle Y_n, \tilde{\mu}(Y_n)\rangle + 2\langle Y_n, \sigma(Y_n)\Delta W_n \rangle + 2\langle Y_n, M_n\rangle
\\& +
2h \langle \tilde{\mu}(Y_n),\sigma (Y_n)\Delta W_n \rangle + 2h \langle \tilde{\mu}(Y_n),M_n \rangle
+2 \langle \sigma(Y_n)\Delta W_n, M_n\rangle
\\ \leq&\|Y_n\|^2 + 3h^2\|\mu(Y_{n})\|^2 + 3\|\sigma(Y_n)\|^2\|\Delta W_n\|^2 +3\|M_n\|^2 \\&  + \frac{2h}{1+h \|\mu(Y_{n})\|}\langle Y_n, \mu(Y_n)\rangle + 2\langle Y_n, \sigma(Y_n)\Delta W_n \rangle + 2\langle Y_n, M_n\rangle
\end{split}
\end{equation}
on $\Omega$ and $ 0\leq n \leq N-1$.
Here we denote
\begin{equation}\label{m_n}
M_n = \frac{1}{2} \sum_{j_1,j_2=1}^{m}L^{j_1}\sigma_{j_2}(Y_{n})(\Delta W^{j_1}_n \Delta W^{j_2}_n-\delta_{j_1,j_2}h).
\end{equation}
Additionally, the global Lipschitz continuity of $\sigma, L^{j_1}\sigma_{j_2}$ implies that for $\|x\|\geq 1$
\begin{equation}\label{es1}
\|\sigma(x)\|^2 \leq (\|\sigma(x)-\sigma(0)\|+\|\sigma(0)\|)^2 \leq (K\|x\|+\|\sigma(0)\|)^2 \leq (K+\|\sigma(0)\|)^2\|x\|^2,
\end{equation}
and
\begin{equation} \label{es2}
\begin{split}
\|L^{j_1}\sigma_{j_2}(x)\|^2 &\leq (\|L^{j_1}\sigma_{j_2}(x)-L^{j_1}\sigma_{j_2}(0)\|+\|L^{j_1}\sigma_{j_2}(0)\|)^2
\\ &\leq (K\|x\|+\|L^{j_1}\sigma_{j_2}(0)\|)^2
\\ &\leq (K+\|L^{j_1}\sigma_{j_2}(0)\|)^2\|x\|^2,
\end{split}
\end{equation}
and the globally one-sided Lipschitz continuity of $\mu$ gives that
\begin{equation} \label{es3}
\langle x, \mu(x) \rangle = \langle x, \mu(x)-\mu(0)\rangle+\langle x,
\mu(0) \rangle \leq K\|x\|^2+ \|x\|\cdot\|\mu(0)\| \leq (K+\|\mu(0)\|)\|x\|^2.
\end{equation}
Furthermore, the polynomial growth bound on $\mu'$ implies that
\begin{equation} \label{es4}
\begin{split}
\|\mu(x)\|^2 &\leq (\|\mu(x)-\mu(0)\|+\|\mu(0)\|)^2 \leq \left(K(1+\|x\|^c)\|x\|+\|\mu(0)\|\right)^2
\\ &\leq \left(2K\|x\|^{(c+1)}+\|\mu(0)\|\right)^2
\\ &\leq \left(2K+\|\mu(0)\|\right)^2\|x\|^{2(c+1)}
\\ &\leq N\left(2K+\|\mu(0)\|\right)^2\|x\|^2
\end{split}
\end{equation}
on $ 1\leq \|x\| \leq N^{\frac{1}{2c}}$. Combining (\ref{es1})-(\ref{es4}), we get from (\ref{es0}) that
\begin{equation}\label{es5}
\begin{split}
\|Y_{n+1}\|^2 \leq&\|Y_n\|^2 + 3hT\left(2K+\|\mu(0)\|\right)^2\|Y_{n}\|^2 + 3(K+\|\sigma(0)\|)^2\|Y_n\|^2\|\Delta W_n\|^2
\\&  +3\|M_n\|^2  + 2h(K+\|\mu(0)\|)\|Y_n\|^2 + 2\langle Y_n, \sigma(Y_n)\Delta W_n \rangle + 2\langle Y_n, M_n\rangle
\end{split}
\end{equation}
on $\{\omega \in \Omega| 1 \leq \|Y_n(\omega)\| \leq N^{\frac{1}{2c}}\}.$ Since $\omega \in \Omega_{n+1}$ implies $\|\Delta W_n\|\leq 1$, on $\omega \in \Omega_{n+1} \cap \{\omega \in \Omega | 1 \leq \|Y_n(\omega)\| \leq N^{\frac{1}{2c}}\}$, we derive from (\ref{m_n}) and (\ref{es2}) that
\begin{equation}\label{es6}
\begin{split}
\|M_n\|^2 \leq& \frac{m^2}{4}\sum_{j_1,j_2=1}^{m}\|L^{j_1}\sigma_{j_2}(Y_n)\|^2|\Delta W^{j_1}_n \Delta W^{j_2}_n-\delta_{j_1,j_2}h|^2
\\ \leq & \frac{m^2}{2}\sum_{j_1,j_2=1}^{m}\left(K+\|L^{j_1}\sigma_{j_2}(0)\|\right)^2\|Y_n\|^2 \cdot\left(|\Delta W^{j_1}_n \Delta W^{j_2}_n|^2+\delta_{j_1,j_2}h^2\right)
\\ \leq & \frac{m^3}{2}\Big(K+\max_{1\leq j_1,j_2\leq m}\|L^{j_1}\sigma_{j_2}(0)\|\Big)^2\left[\|\Delta W_n\|^2 + h^2\right]\|Y_n\|^2,
\end{split}
\end{equation}
where the fact that $|\Delta W_n^{j_2}|^2 \leq \|\Delta W_n\|^2 \leq 1 $ was used. Inserting (\ref{es6}) into (\ref{es5}) shows that
\begin{equation} \label{yn2}
\begin{split}
&\|Y_{n+1}\|^2
\\ \leq& \|Y_n\|^2 \Big[1+ 3 \frac{T^2}{N}(2K + \|\mu(0)\|)^2 + 3(K+\|\sigma(0)\|)^2\|\Delta W_n\|^2+2\frac{T}{N}(K+\|\mu(0)\|)\Big]
\\ & + \frac{3m^3}{2}\Big(K+\max_{1\leq j_1,j_2\leq m}\|L^{j_1}\sigma_{j_2}(0)\|\Big)^2\, \Big(\|\Delta W_n\|^2 + h^2\Big)\|Y_n\|^2
+ 2\langle Y_n, \sigma(Y_n)\Delta W_n \rangle + 2\langle Y_n, M_n\rangle
\\ \leq&\|Y_n\|^2\left[1+ \frac{2}{N}\left(\frac{3}{2}(2TK+T\|\mu(0)\|)^2 + T(K+\|\mu(0)\|) +\frac{3m^3T^2}{4}\Big(K+\max_{1\leq j_1,j_2\leq m}\|L^{j_1}\sigma_{j_2}(0)\|\Big)^2\right) \right.
\\ &+ \left.2\left(\frac{3}{2}(K+\|\sigma(0)\|)^2 +\frac{3m^3}{4}\Big(K+\max_{1\leq j_1,j_2\leq m}\|L^{j_1}\sigma_{j_2}(0)\|\Big)^2\right)\|\Delta W_n\|^2 + 2 \alpha_n\right]
\\ \leq& \|Y_n\|^2 \exp\Big[\frac{2\lambda}{N} + 2\lambda \|\Delta W_n\|^2 + 2\alpha_n \Big]
\end{split}
\end{equation}
on $\Omega_{n+1} \cap \{\omega \in \Omega | 1 \leq \|Y_n(\omega)\| \leq N^{\frac{1}{2c}}\}$.
Now combining (\ref{yn1}) and (\ref{yn2}), and using mathematical induction as used in the proof of Lemma 2.1 in \cite{HJK10} finish the proof.
$\square$

\begin{lem} \label{lemma0}
For all $p \geq 1$
\begin{equation}\label{d0}
\sup_{\begin{subarray}{ll}N\in\mathbb{N}\\ N \geq 4\lambda pT \end{subarray}}\mathbb{E}\Big[ \exp\Big(p\lambda\sum_{k=0}^{N-1}\|\Delta W_k\|^2\Big)\Big] < \infty.
\end{equation}
\end{lem}
{\it Proof.} This result is identical to Lemma 2.3 in \cite{HJK10} with only different $\lambda$. $\square$

\begin{lem} \label{lemma1}
Let $\alpha_n$ be given by (\ref{alpha}). Then for all $p \geq 1$
\begin{equation} \label{d2}
\sup_{z\in\{-1,1\}}\sup_{N\in\mathbb{N}}\mathbb{E}\Big[\sup_{0\leq n\leq N} \exp\Big(pz\sum_{k=0}^{n-1}\alpha_k\Big)\Big] < \infty.
\end{equation}
\end{lem}
{\it Proof.} We set
$$
\alpha_k^1 = 1_{\{\|Y_k\| \geq 1\}}\Big\langle\frac{Y_k}{\|Y_k\|}, \frac{\sigma(Y_k)}{\|Y_k\|}\Delta W_k\Big\rangle,
$$
$$\alpha_k^2 = 1_{\{\|Y_k\| \geq 1\}}\Big\langle\frac{Y_k}{\|Y_k\|}, \sum_{j_1,j_2=1}^m\frac{L^{j_1}\sigma_{j_2}(Y_k)}{2\|Y_k\|}\left(\Delta W^{j_1}_k \Delta W^{j_2}_k-\delta_{j_1,j_2}h\right)\Big\rangle.
$$
Then we have $\alpha_k = \alpha_k^1 + \alpha_k^2$. Hence H\"{o}lder's inequality shows that
\begin{equation}\label{alphab1}
\Big\|\sup_{0\leq n\leq N}\exp \Big( z \sum_{k=0}^{n-1} \alpha_k \Big)\Big\|_{L^p(\Omega;\mathbb{R})} \leq \Big\|\sup_{0\leq n\leq N} \exp \Big( z \sum_{k=0}^{n-1} \alpha_k^1 \Big)\Big\|_{L^{2p}(\Omega;\mathbb{R})} \cdot \Big\|\sup_{0\leq n\leq N} \exp \Big( z \sum_{k=0}^{n-1} \alpha_k^2 \Big)\Big\|_{L^{2p}(\Omega;\mathbb{R})}.
\end{equation}
Note that Lemma 2.4 in \cite{HJK10} has proved that
\begin{equation}
\sup_{z\in\{-1,1\}}\Big\|\sup_{0\leq n\leq N} \exp \Big( z \sum_{k=0}^{n-1} \alpha_k^1 \Big)\Big\|_{L^{2p}(\Omega;\mathbb{R})} < \infty.
\end{equation}
Consequently it remains to prove the boundedness of the second term on the right-hand side of (\ref{alphab1}). One can easily verify that the discrete stochastic process $z\sum_{k=0}^{n-1} \alpha^2_k, n \in \{0,1,...,N\}$ is an $\{\mathcal {F}_{t_n}: 0\leq n \leq N\}$-martingale for every $ z \in \{-1,1\}$. Since the exponential function is convex, the discrete stochastic process $\exp\left(z\sum_{k=0}^{n-1} \alpha^2_k\right), n \in \{0,1,...,N\}$ is a positive  $\{\mathcal {F}_{t_n}: 0\leq n \leq N\}$-submartingale for every $z \in \{-1,1\}$. Therefore, Doob's maximal inequality gives that
\begin{equation}\label{dr1}
\Big\|\sup_{0\leq n\leq N} \exp \Big( z \sum_{k=0}^{n-1} \alpha^2_k \Big)\Big\|_{L^{2p}(\Omega;\mathbb{R})} \leq \frac{2p}{2p-1} \Big\|\exp \Big( z \sum_{k=0}^{N-1} \alpha^2_k \Big)\Big\|_{L^{2p}(\Omega;\mathbb{R})}.
\end{equation}
The Cauchy-Schwarz inequality, (\ref{es2}) and Lemma \ref{lemma0} give that
\begin{eqnarray*}
%\begin{split}
&&\Big\|\exp \Big( z \sum_{k=0}^{N-1} \alpha_k^2 \Big)\Big\|_{L^{2p}(\Omega;\mathbb{R})}
\\ &\leq&
\Big\|\exp \Big(  \sum_{k=0}^{N-1} 1_{\{\|Y_k\| \geq 1\}} \sum_{j_1,j_2=1}^m\frac{\|L^{j_1}\sigma_{j_2}(Y_k)\|}{2\|Y_k\|}\Big|\Delta W^{j_1}_k \Delta W^{j_2}_k-\delta_{j_1,j_2}h\Big| \Big)\Big\|_{L^{2p}(\Omega;\mathbb{R})}
\\ &\leq&
\Big\|\exp \Big( \sum_{k=0}^{N-1} (\frac{K}{2}+\frac{1}{2}\max_{1\leq j_1,j_2\leq m}\|L^{j_1}\sigma_{j_2}(0)\|) \sum_{j_1,j_2=1}^m\Big(\Big|\Delta W^{j_1}_k \Delta W^{j_2}_k\Big|+\Big|\delta_{j_1,j_2}h\Big|\Big) \Big)\Big\|_{L^{2p}(\Omega;\mathbb{R})}
\\ &\leq&
\Big\|\exp \Big(  \sum_{k=0}^{N-1} (\frac{K}{2}+\frac{1}{2}\max_{1\leq j_1,j_2\leq m}\|L^{j_1}\sigma_{j_2}(0)\|)\Big( m \|\Delta W_k\|^2 +mh \Big) \Big)\Big\|_{L^{2p}(\Omega;\mathbb{R})}
\\ &\leq&
e^{\lambda}\Big\|\exp \Big( \lambda \sum_{k=0}^{N-1}  \|\Delta W_k\|^2 \Big)\Big\|_{L^{2p}(\Omega;\mathbb{R})}
< \infty.
%\end{split}
\end{eqnarray*}
This together with (\ref{dr1}) completes the proof.
$\square$

\begin{lem} \label{lemma2}
Let $D_n$ be given by (\ref{Dn}).  Then for all $N \geq 8\lambda pT$ and $p \geq 1$
\begin{equation} \label{d03}
\sup_{N\in\mathbb{N}}\mathbb{E}\left[\sup_{0\leq n\leq N} |D_n|^p\right] < \infty.
\end{equation}
\end{lem}
{\it Proof.} %From (\ref{Dn}) we have
%\begin{equation}\nonumber
%D_n = e^{\lambda}(\lambda + \|\xi\|)\cdot \exp\left(\lambda\sup_{0\leq u\leq n}\sum_{k=u}^{n-1}\|\Delta W_k\|^2 \right)\cdot \exp\left(\sup_{0\leq u\leq n}\left[\sum_{k=0}^{n-1}\alpha_k- \sum_{k=0}^{u-1}\alpha_k \right]\right).
%\end{equation}
%Thus
%\begin{equation}\nonumber
%\begin{split}
%\sup_{0\leq n\leq N}|D_n|^p \leq& e^{p\lambda}(\lambda + \|\xi\|)^p\cdot \exp\left(p \lambda\sum_{k=0}^{N-1}\|\Delta W_k\|^2 \right)
%\\ &\cdot \sup_{0\leq n\leq N}\left[ \exp\left(p\sum_{k=0}^{n-1}\alpha_k\right)\right] \cdot \sup_{0\leq u\leq N}\left[ \exp\left(-p\sum_{k=0}^{u-1}\alpha_k\right)\right].
%\end{split}
%\end{equation}
Note that $D_n$ here takes the same form as $D_n^N$ in \cite{HJK10}, with only different $\alpha_n$. With Lemma \ref{lemma0} and Lemma \ref{lemma1} at hand, one can follow the proof of Lemma 2.5 in \cite{HJK10} to derive the desired result. $\square$

\begin{lem} \label{lemma3}
Let $\Omega_N$ be given by (\ref{Omega}) with $n = N$.  Then for all $p \geq 1$
\begin{equation} \label{d3}
\sup_{N\in\mathbb{N}}\left(N^p\mathbb{P}\left[(\Omega_N)^c\right]\right) < \infty.
\end{equation}
\end{lem}
{\it Proof.} The proof is identical to the proof of Lemma 2.6 in \cite{HJK10}. $\square$
%Since $\Omega_n, D_n$ here share the same form as $\Omega_n, D_n$ in \cite{HJK10},
%the proof of

Before establishing the main result of this section, we also need two Burkholder-Davis-Gundy
type inequalities.

\begin{lem}\label{BDG-Continuous}
Let $k \in \mathbb{N}$ and let $Z: [0,T] \times \Omega \rightarrow \mathbb{R}^{k \times m}$ be a predictable stochastic process satisfying $\mathbb{P}(\int_0^T \|Z_s\|^2 ds < \infty) = 1$. Then for all $t \in [0,T]$ and all $p \geq 2$
\begin{equation}\label{BDGc}
\Big\|\sup_{s\in[0,t]}\Big\|\int_0^s Z_u d W_u\Big\|\Big\|_{L^p(\Omega;\mathbb{R})} \leq p \Big(\int_0^t \sum_{i=1}^m \|Z_s e_i\|_{L^p(\Omega;\mathbb{R}^k)}^2 ds\Big)^{\frac{1}{2}}.
\end{equation}
Here the vectors $e_1 = (1,0,...,0)^T \in \mathbb{R}^m$, $e_2 = (0,1,...,0)^T \in \mathbb{R}^m,...,$ $e_m = (0,0,...,1)^T \in \mathbb{R}^m$ are orthogonal basis of vector space $\mathbb{R}^m$.
\end{lem}
{\it Proof.} Combining Doob's maximal inequality and Lemma 7.7 of Da Prato, G., and Zabczyk \cite{DZ92} gives the desired assertion. $\square$

The following is a discrete version of the Burkholder-Davis-Gundy
type inequality (\ref{BDGc}).
\begin{lem}\label{BDG-Discrete}
Let $k \in \mathbb{N}$ and let $Z_l: \Omega \rightarrow \mathbb{R}^{k \times m}, l \in \{0,1,...,N-1\}$ be a family of mappings such that  $Z_l$ is $\mathcal{F}_{\frac{lT}{N}}/\mathcal{B}(\mathbb{R}^{k\times m})$-measurable. Then for all $0\leq n\leq N$ and $p \geq 2$
\begin{equation}\label{BDGd}
\Big\|\sup_{0\leq j\leq n}\Big\|\sum_{l=0}^{j-1} Z_l \Delta W_l\Big\|\Big\|_{L^p(\Omega;\mathbb{R})} \leq p \Big(\sum_{l=0}^{n-1} \sum_{i=1}^m \|Z_l e_i\|_{L^p(\Omega;\mathbb{R}^k)}^2 \frac{T}{N}\Big)^{\frac{1}{2}}.
\end{equation}
\end{lem}

\begin{thm} \label{lemma4}
Let $Y_n$  be given by (\ref{MilsteinC}).  Then for all $p \in [1,\infty)$
\begin{equation} \label{d4}
\sup_{N\in\mathbb{N}}\left[\sup_{0\leq n\leq N}\mathbb{E} \|Y_n\|^p\right] < \infty.
\end{equation}
\end{thm}
{\it Proof.} From (\ref{MilsteinC}) we have
\begin{equation} \label{mb1}
\begin{split}
&\|Y_n\|_{L^p(\Omega;\mathbb{R}^d)} \leq \|\xi\|_{L^p(\Omega;\mathbb{R}^d)} \\ &+
\Big\|\sum_{k=0}^{n-1} h\tilde{\mu}(Y_{k})\Big\|_{L^p(\Omega;\mathbb{R}^d)}
+ \Big\|\sum_{k=0}^{n-1} \sigma (Y_k)\Delta W_k\Big\|_{L^p(\Omega;\mathbb{R}^d)}
 +\Big\|\sum_{k=0}^{n-1}M_k\Big\|_{L^p(\Omega;\mathbb{R}^d)},
\end{split}
\end{equation}
where the notation $M_k$ comes from (\ref{m_n}). Using \eqref{OLC3}, the triangle inequality and the Burkholder-Davis-Gundy
type inequality in Lemma \ref{BDG-Discrete} we have
\begin{equation} \label{mb2}
\begin{split}
\Big\|\sum_{k=0}^{n-1} \sigma (Y_k)\Delta W_k\Big\|_{L^p(\Omega;\mathbb{R}^d)} \leq &
\Big\|\sum_{k=0}^{n-1} \left[\sigma (Y_k)-\sigma (0)\right]\Delta W_k\Big\|_{L^p(\Omega;\mathbb{R}^d)} + \Big\|\sum_{k=0}^{n-1} \sigma (0) \Delta W_k\Big\|_{L^p(\Omega;\mathbb{R}^d)}
\\ \leq &  p\Big(\sum_{k=0}^{n-1}\sum_{i=1}^{m} \left\|\sigma_i(Y_k)-\sigma_i(0)\right\|^2_{L^p(\Omega;\mathbb{R}^d)}h\Big)^{1/2}
+ p\Big(\frac{nT}{N}\sum_{i=1}^{m} \left\|\sigma_i(0)\right\|^2\Big)^{1/2}
\\ \leq &  p \Big(mK^2h \sum_{k=0}^{n-1} \left\|Y_k\right\|^2_{L^p(\Omega;\mathbb{R}^d)}\Big)^{1/2}
+ p\sqrt{Tm} \|\sigma(0)\|.
\end{split}
\end{equation}
For the fourth term on the right-hand side of (\ref{mb1}), the estimate in the second inequality of (\ref{es2})  and the independence of $Y_k$ and $\Delta W_k$ imply that
\begin{equation} \label{mb3}
\begin{split}
&\Big\|\sum_{k=0}^{n-1}M_k\Big\|_{L^p(\Omega;\mathbb{R}^d)} =\frac{1}{2}\Big\|\sum_{k=0}^{n-1}\Big(\sum_{j_1,j_2=1}^{m}L^{j_1}\sigma_{j_2}(Y_k)
\left(\Delta W^{j_1}_k \Delta W^{j_2}_k-\delta_{j_1,j_2}h\right) \Big)\Big\|_{L^p(\Omega;\mathbb{R}^d)}
\\ \leq &  \frac{1}{2}\sum_{k=0}^{n-1}\sum_{j_1,j_2=1}^{m}\left[ \left(K\left\|Y_k\right\|_{L^p(\Omega;\mathbb{R}^d)}+\|L^{j_1}\sigma_{j_2}(0)\|\right)
\cdot\left\|\Delta W^{j_1}_k \Delta W^{j_2}_k-\delta_{j_1,j_2}h\right\|_{L^p(\Omega;\mathbb{R})}
\right].
\end{split}
\end{equation}
Using the Burkholder-Davis-Gundy
type inequality \eqref{BDGd} and mutual independence of $\Delta W_k$ gives
\begin{equation} \label{mb4}
\begin{split}
&\frac{1}{2}\sum_{j_1,j_2=1}^{m} \|\Delta W^{j_1}_k \Delta W^{j_2}_k-\delta_{j_1,j_2}h\|_{L^p(\Omega;\mathbb{R})}
\leq  \frac{1}{2}\sum_{j_1,j_2=1}^{m} \|\Delta W^{j_1}_k \Delta W^{j_2}_k\|_{L^p(\Omega;\mathbb{R})} + \frac{1}{2}mh
\\ = & \frac{1}{2} \sum_{\begin{subarray}{ll}j_1,j_2=1\\j_1 \neq j_2 \end{subarray}}^{m}\|\Delta W_k^{j_1}\|_{L^p(\Omega;\mathbb{R})}\cdot \|\Delta W_k^{j_2}\|_{L^p(\Omega;\mathbb{R})} + \frac{1}{2}\sum_{j=1}^{m}\|\Delta W_k^{j}\|^2_{L^{2p}(\Omega;\mathbb{R})} + \frac{1}{2}mh
\\ \leq & c_{p,m}h,
\end{split}
\end{equation}
where $c_{p,m}=\frac{m^2-m}{2}p^2 + 2mp^2 + \frac{m}{2}$. Inserting (\ref{mb4}) into (\ref{mb3}) we obtain that
\begin{equation} \label{mb5}
\begin{split}
\Big\|\sum_{k=0}^{n-1}M_k\Big\|_{L^p(\Omega;\mathbb{R}^d)}\leq
  Kc_{p,m}h\sum_{k=0}^{n-1}\|Y_k\|_{L^p(\Omega;\mathbb{R}^d)} + Tc_{p,m}\sup_{1\leq j_1,j_2\leq m}\|L^{j_1}\sigma_{j_2}(0)\|.
\end{split}
\end{equation}
Combining (\ref{mb1}),(\ref{mb2}) and (\ref{mb5}), and using $\|\tilde{\mu}(Y_k)\| \leq 1$ give
\begin{equation}
\begin{split}
\|Y_n\|_{L^p(\Omega;\mathbb{R}^d)}
\leq &\|\xi\|_{L^p(\Omega;\mathbb{R}^d)} + N + p\Big(mK^2h \sum_{k=0}^{n-1}\|Y_k\|^2_{L^p(\Omega;\mathbb{R}^d)}\Big)^{1/2} + p\sqrt{Tm} \|\sigma(0)\|
\\& + Kc_{p,m}h\sum_{k=0}^{n-1}\|Y_k\|_{L^p(\Omega;\mathbb{R}^d)} + Tc_{p,m}\sup_{1\leq j_1,j_2\leq m}\|L^{j_1}\sigma_{j_2}(0)\|.
\end{split}
\end{equation}
Thus taking square of both sides shows that
\begin{equation}
\begin{split}
\|Y_n\|_{L^p(\Omega;\mathbb{R}^d)}^2 \leq& 3\Big(\|\xi\|_{L^p(\Omega;\mathbb{R}^d)} + N + p\sqrt{Tm} \|\sigma(0)\|+ Tc_{p,m}\sup_{1\leq j_1,j_2\leq m}\|L^{j_1}\sigma_{j_2}(0)\|\Big)^2
\\&+ 3mK^2p^2h \sum_{k=0}^{n-1}\|Y_k\|_{L^p(\Omega;\mathbb{R}^d)}^2
+3K^2c_{m,p}^2Th\sum_{k=0}^{n-1}\|Y_k\|^2_{L^p(\Omega;\mathbb{R}^d)}.
\end{split}
\end{equation}
In the next step Gronwall's lemma gives that
\begin{equation} \label{le3}
\begin{split}
\sup_{0\leq n\leq N}\|Y_n\|_{L^p(\Omega;\mathbb{R}^d)}
\leq \sqrt{3}\exp(C_1)\left(\|\xi\|_{L^p(\Omega;\mathbb{R}^d)} + N + C_2\right),
\end{split}
\end{equation}
where
$
C_1 = \frac{1}{2}(3mK^2p^2+3K^2c_{p,m}^2T)T, \, C_2 = p\sqrt{Tm} \|\sigma(0)\|+ Tc_{p,m}\sup_{1\leq j_1,j_2\leq m}\|L^{j_1}\sigma_{j_2}(0)\|.
$
Due to the $N$ on the right-hand side of (\ref{le3}), (\ref{le3}) does not complete the prove. However, exploiting (\ref{le3}) in an
appropriate bootstrap argument will enable us to establish (\ref{d4}). First, H\"{o}lder's inequality,  Lemma \ref{lemma3} and the estimate (\ref{le3}) show that
\begin{equation}\label{le4}
\begin{split}
&\sup_{N\in\mathbb{N}}\sup_{0\leq n \leq N}\|
1_{(\Omega_n)^c}Y_n\|_{L^p(\Omega;\mathbb{R}^d)}
\\ \leq &
\sup_{N\in\mathbb{N}}\sup_{0\leq n \leq N} \left(\|
1_{(\Omega_n)^c}\|_{L^{2p}(\Omega;\mathbb{R}^d)}\|Y_n\|_{L^{2p}(\Omega;\mathbb{R}^d)}\right)
\\ \leq &
\Big(\sup_{N\in\mathbb{N}}\left(N\cdot \|
1_{(\Omega_N)^c}\|_{L^{2p}(\Omega;\mathbb{R}^d)}\right)\Big)\Big(\sup_{N\in\mathbb{N}}\sup_{0\leq n \leq N} \left(N^{-1}\cdot\|Y_n\|_{L^{2p}(\Omega;\mathbb{R}^d)}\right)\Big)
\\ \leq &
\Big(\sup_{N\in\mathbb{N}}N^{2p}\cdot \mathbb{P}\left[(\Omega_N)^c\right]\Big)^{\frac{1}{2p}}\Big(\sup_{N\in\mathbb{N}}\sup_{0\leq n \leq N} \left(N^{-1}\cdot\|Y_n\|_{L^{2p}(\Omega;\mathbb{R}^d)}\right)\Big)
\\ < &\infty.
\end{split}
\end{equation}
In addition, Lemma \ref{DomLem} and Lemma \ref{lemma2} imply that
\begin{equation}\label{le5}
\sup_{N\in\mathbb{N}}\sup_{0\leq n \leq N}\|
1_{\Omega_n}Y_n\|_{L^p(\Omega;\mathbb{R}^d)} \leq \sup_{N\in\mathbb{N}}\sup_{0\leq n \leq N}\|
D_n\|_{L^p(\Omega;\mathbb{R}^d)} < \infty.
\end{equation}
Combining (\ref{le4}) and (\ref{le5}) finally completes the proof.
$\square$

\section{Strong convergence order of the tamed Milstein method}

In order to recover the strong convergence order one for the tamed Milstein method, we additionally need the following assumptions. Throughout this section $C_{p,T}$ is a generic constant that might vary from one place to another and depends on $\mu,\sigma$, the initial data $\xi$, and the interval of integration $[0,T]$, but is independent of the discretisation parameter.
\begin{ass} \label{SCC}
Assume that $\mu(x)$ and $\sigma_{i}(x)$ are two times continuous differentiable and there exist positive constants $K,q \geq 1$, such that $\forall
x \in \mathbb{R}^d, i=1,...,d$
\begin{eqnarray}
\|\mu''(x)\|_{L^{(2)}(\mathbb{R}^d;\mathbb{R}^d)} &\leq& K(1+\|x\|^q),\label{OLC6}\\
\|\sigma''_{i}(x)\|_{L^{(2)}(\mathbb{R}^d;\mathbb{R}^d)} &\leq& K.    \label{OLC7}
\end{eqnarray}
\end{ass}
Here, for a two times continuous differentiable function $f:\mathbb{R}^d\rightarrow \mathbb{R}^d$ we use the notation $\|f''(x)\|_{L^{(2)}(\mathbb{R}^d;\mathbb{R}^d)} = \sup_{h_1,h_2 \in \mathbb{R}^d, \|h_1\|\leq 1, \|h_2\|\leq 1 }\|f''(x)(h_1,h_2)\|$. The bilinear operator $f''(x): \mathbb{R}^d \times \mathbb{R}^d \rightarrow\mathbb{R}^d$ is defined by (\ref{derivative}).
In the following convergence analysis, we prefer to write the tamed Milstein method in the form of (\ref{TamedMilstein}), rather than \eqref{MilsteinC}. We now introduce appropriate time continuous interpolations of the time
discrete numerical approximations. More accurately, we define the time continuous approximation $\bar{Y}_s$ such that for $s \in [t_n,t_{n+1})$
\begin{equation}  \label{Milstein2}
\begin{split}
\bar{Y}_s := & Y_{n} + (s-t_n)\tilde{\mu}(Y_{n}) + \sigma (Y_{n})(W_s- W_{t_n}) + \sum^m_{j_1,j_2=1}L^{j_1}\sigma_{j_2}(Y_{n})I_{j_1,j_2}^{t_n,s}
\\= &Y_{n} + \int_{t_n}^{s} \tilde{\mu} (Y_n)du + \sum_{i=1}^{m}\int_{t_n}^{s}\sigma_{i} (Y_{n})d W_u^{i} + \sum^m_{i=1} \int^{s}_{t_n} \sum_{j=1}^{m}L^{j}\sigma_{i}(Y_{n})\Delta W_{u}^{j}
d W_{u}^{i},
\end{split}
\end{equation}
where we use the notation
$$
I_{j_1,j_2}^{t_n,s} = \int^{s}_{t_n} \int^{s_2}_{t_n} d W_{s_1}^{j_1}
d W_{s_2}^{j_2}, \quad \Delta W_{s}^{j} := \sum_{n=0}^{\infty}1_{\{t_n\leq s< t_{n+1}\}} (W_{s}^{j}-W_{t_n}^{j}).
$$
It is evident that $\bar{Y}_{t_n} = Y_n, n= 0,1,...,N$, that is,
$\bar{Y}_t$ coincides with the discrete solutions at the grid-points.
We can rewrite $\bar{Y}_t$ as an integral form in the whole interval $[0,T]$
\begin{equation}  \label{MilCon}
\begin{split}
\bar{Y}_t = &Y_{0} + \int_{0}^{t} \tilde{\mu}(Y_{n_s}) ds + \sum_{i=1}^{m}\int_{0}^{t}\Big[\sigma_{i}(Y_{n_s})+\sum_{j=1}^{m}L^{j}\sigma_{i}(Y_{n_s})\Delta W_{s}^{j}\Big]d W_s^{i},
\end{split}
\end{equation}
where $n_s$ is the greatest integer number such that $t_{n_s} \leq s$.
%\begin{equation} \nonumber
%Y_{n_s} := \sum_{n=0}^{\infty}1_{\{t_n\leq s< t_{n+1}\}} Y_n.
%\end{equation}
Combining (\ref{itsde}) and (\ref{MilCon}) gives
\begin{equation}\label{xyd}
\begin{split}
X_t- \bar{Y}_t =&  \int_{0}^{t}\left[\mu(X_s)- \tilde{\mu}(Y_{n_s}) \right]ds+ \sum_{i=1}^{m}\int_{0}^{t}\Big[\sigma_{i}(X_s)-\sigma_{i}(Y_{n_s})-\sum_{j=1}^{m}L^{j}\sigma_{i}(Y_{n_s})\Delta W_{s}^{j}\Big]d W_s^{i}.
\end{split}
\end{equation}

In what follows, we also use deterministic Taylor formula frequently. If a function  $f:\mathbb{R}^d \rightarrow \mathbb{R}^d$ is twice differentiable, the following Taylor's formula holds \cite{AP93}
\begin{equation}\label{taylorformula1}
\begin{split}
&f(\bar{Y}_s)-f(Y_{n_s}) = f'(Y_{n_s})(\bar{Y}_s-Y_{n_s}) + R_1(f),
\end{split}
\end{equation}
where $R_1(f)$ is the remainder term
\begin{equation}\label{R1}
\begin{split}
R_1(f) =& \int_0^1(1-r) f''(Y_{n_s}+r(\bar{Y}_s-Y_{n_s}))(\bar{Y}_s-Y_{n_s},\bar{Y}_s-Y_{n_s}) dr.
\end{split}
\end{equation}
Here for arbitrary $a, h_1,h_2 \in \mathbb{R}^d$ the derivatives have the following expression
\begin{equation}\label{derivative}
f'(a)(h_1) = \sum_{i=1}^{d} \frac{\partial f}{\partial x^{i}}h_1^{i}, \quad
f''(a)(h_1,h_2) = \sum_{i,j=1}^{d} \frac{\partial^2 f}{\partial x^{i}\partial x^{j}}h_1^{i}h_2^{j}.
\end{equation}
Replacing $\bar{Y}_s-Y_{n_s}$ in \eqref{taylorformula1} with \eqref{Milstein2} and rearranging lead to
\begin{equation}\label{taylorformula2}
\begin{split}
f(\bar{Y}_s)-f(Y_{n_s})  = f'(Y_{n_s})\big(\sigma(Y_{n_s})(W_s- W_{t_{n_s}})\big) + \tilde{R}_1(f),
\end{split}
\end{equation}
where
\begin{equation}\label{RT1}
\begin{split}
\tilde{R}_1(f)=f'(Y_{n_s}) \Big((s-t_{n_s}) \tilde{\mu}(Y_{n_s})
+\sum^m_{j_1,j_2=1}L^{j_1}\sigma_{j_2}(Y_{n_s})I_{j_1,j_2}^{t_{n_s},s}\Big) +R_1(f).
\end{split}
\end{equation}

By the definitions (\ref{Lj}) and (\ref{derivative}), it can be readily checked that
\begin{equation}\label{eq:sigmac}
\sigma_i'(x) \big(\sigma_j(x) \big)= L^{j}\sigma_i(x).
\end{equation}
Therefore replacing $f$ in \eqref{taylorformula2} by $\sigma_{i}$ and taking  \eqref{eq:sigmac} into account show that
\begin{equation} \label{eq:sigmaTaylor}
\begin{split}
\tilde{R}_1(\sigma_{i}) = \sigma_{i}(\bar{Y}_s)-\sigma_{i}(Y_{n_s})-\sum_{j=1}^{m}L^{j}\sigma_{i}(Y_{n_s})\Delta W_{s}^{j}.
\end{split}
\end{equation}

\begin{thm} \label{conthm}
Let conditions in Assumptions \ref{LCMC} and \ref{SCC} and (\ref{cc}) be fulfilled.  Then there exists a family of real numbers $C_{p,T} \geq 1$, $p\geq 1$ such that
\begin{equation} \label{d5}
\left(\mathbb{E} \Big[\sup_{t\in[0,T]}\|X_t -\bar{Y}_t\|^p\Big]\right)^{1/p} \leq C_{p,T}\cdot h, \quad h \in (0, 1].
\end{equation}
\end{thm}
The following five lemmas are used in the proof of Theorem \ref{conthm}.
\begin{lem} \label{mblem0}
Let conditions in Assumptions \ref{LCMC} be fulfilled.  Then for all $p\geq 1, 1\leq j_1,j_2 \leq m$ the following estimates hold
\begin{equation}
\sup_{N \in \mathbb{N}}\sup_{0\leq n \leq N} \left[ \mathbb{E}\|\mu(Y_n)\|^{p} \bigvee \mathbb{E} \|\mu'(Y_n)\|^{p} \bigvee \mathbb{E} \|\sigma(Y_n)\|^{p} \bigvee \mathbb{E} \|L^{j_1} \sigma_{j_2}(Y_n)\|^p\right] < \infty.
\end{equation}
\end{lem}
{\it Proof.} It immediately follows from Theorem \ref{lemma4} by considering (\ref{OLC5}) and (\ref{es1})-(\ref{es4}). $\square$

\begin{lem} \label{mblem}
Let conditions in Theorem \ref{conthm} be fulfilled. Then for all $p \geq 1$
\begin{equation} \label{d7}
\sup_{t\in[0,T]} \Big[\|X_t \|_{L^p(\Omega;\mathbb{R}^d)} \bigvee \|\mu(X_t) \|_{L^p(\Omega;\mathbb{R}^d)} \bigvee  \|\sigma(X_t) \|_{L^p(\Omega;\mathbb{R}^{d\times m})} \bigvee  \|\bar{Y}_t \|_{L^p(\Omega;\mathbb{R}^d)}\Big] < \infty.
\end{equation}
\end{lem}
{\it Proof}. Conditions (\ref{OLC1})-(\ref{OLC3}) ensure that the exact solution $X_t$ satisfies $\sup_{t\in[0,T]}\mathbb{E} \|X_t \|^{p} < \infty$ (See, e.g., \cite[Theorem 4.1]{Mao97}). Polynomial growth condition on $\mu$ as (\ref{es4}) and linear growth condition on $\sigma$ give the second and the third estimates. Taking the definition (\ref{Milstein2}) and Lemma \ref{mblem0} into account, we can easily obtain the last estimate. $\square$
\begin{lem} \label{yd}
Let conditions in Theorem \ref{conthm} be fulfilled.  Then there exists a family of constants $C_{p,T} \geq 1$ such that for $p \geq 1$
\begin{equation} \label{d6}
\sup_{t\in[0,T]}\|\bar{Y}_t - Y_{n_t}\|_{L^{p}(\Omega;\mathbb{R}^d)} \bigvee  \sup_{t\in[0,T]}\|\mu(X_t) - \mu(X_{t_{n_t}})\|_{L^{p}(\Omega;\mathbb{R}^d)} \leq C_{p,T}\, h^{\frac{1}{2}}.
\end{equation}
\end{lem}
{\it Proof}. Let $n_t$ be the greatest integer number such that $t_{n_t} \leq t$. From (\ref{Milstein2}) we have
\begin{equation}\label{ty}
\bar{Y}_t - Y_{n_t} = (t-t_{n_t}) \tilde{\mu}(Y_{n_t}) + \sigma (Y_{n_t})(W_t- W_{t_{n_t}}) + \sum^m_{j_1,j_2=1}L^{j_1}\sigma_{j_2}(Y_{n_t})I_{j_1,j_2}^{t_{n_t},t}.
\end{equation}
Following the same line as estimating (\ref{mb1}), one can readily derive the first estimate. For the second estimate, we use \eqref{OLC5} and H\"{o}lder's inequality to obtain
\begin{equation}\label{tx}
\begin{split}
\|\mu(X_t) - \mu(X_{t_{n_t}})\|_{L^{p}(\Omega;\mathbb{R}^d)}  \leq&  K \big\| (1+ \|X_t\|^c + \|X_{t_{n_t}}\|^c)\cdot \|X_t - X_{t_{n_t}}\|  \big\|_{L^{p}(\Omega;\mathbb{R})}
\\ \leq &
K \big\| (1+ \|X_t\|^c + \|X_{t_{n_t}}\|^c)\big\|_{L^{2p}(\Omega;\mathbb{R})} \cdot \big\|X_t - X_{t_{n_t}}\big\|_{L^{2p}(\Omega;\mathbb{R}^d)},
\end{split}
\end{equation}
where
\begin{equation}\label{eq:dx}
X_t - X_{t_{n_t}} = \int_{t_{n_t}}^t \mu(X_s) ds + \int_{t_{n_t}}^t \sigma(X_s) dW_s.
\end{equation}
Taking Lemma \ref{mblem} into account and using the same argument as in the previous section, one can obtain $\big\|X_t - X_{t_{n_t}}\big\|_{L^{2p}(\Omega;\mathbb{R}^d)} \leq C_{p,T} h^{\frac{1}{2}}$, and then complete the proof easily by considering \eqref{tx} and Lemma \ref{mblem}. $\square$
\begin{lem} \label{Rlem}
Let $p \geq 1$ and let conditions in Theorem \ref{conthm} be fulfilled. Then for $i =1,2,..., m$
\begin{equation} \label{R_estimate}
\|\tilde{R}_1(\mu) \|_{L^p(\Omega;\mathbb{R}^d)} \vee \|\tilde{R}_1(\sigma_{i})\|_{L^p(\Omega;\mathbb{R}^d)} \leq C_{p,T} h.
\end{equation}
\end{lem}
{\it Proof}.
To estimate $\|\tilde{R}_1(\mu)\|_{L^{p}(\Omega;\mathbb{R}^d)}$, we need to estimate $\|R_1(\mu)\|_{L^{p}(\Omega;\mathbb{R}^d)}$ first. Due to Theorem \ref{lemma4} and Lemma \ref{mblem} we can find some suitable constant $C_{p,T}$ such that
\begin{equation}\label{R1mu}
\begin{split}
\|R_1(\mu)\|_{L^{p}(\Omega;\mathbb{R}^d)}
 \leq &  \int_0^1 (1-r) \Big\|\mu^{''}(Y_{n_s}+r(\bar{Y}_s-Y_{n_s}))\, \left(\bar{Y}_s-Y_{n_s},\bar{Y}_s-Y_{n_s}\right)\Big\|_{L^{p}(\Omega;\mathbb{R})} dr
\\ \leq &  \int_0^1 \left\|\|\mu^{''}(Y_{n_s}+r(\bar{Y}_s-Y_{n_s}))\|_{L^{(2)}(\mathbb{R}^d;\mathbb{R}^d)} \cdot\|\bar{Y}_s-Y_{n_s}\|^2\right\|_{L^{p}(\Omega;\mathbb{R})}dr
\\ \leq &   K\left\|\left(1+\|Y_{n_s}\|^q+\|\bar{Y}_s\|^q\right) \cdot\|\bar{Y}_s-Y_{n_s}\|^2\right\|_{L^{p}(\Omega;\mathbb{R})}
\\ \leq &   K\left\|1+\|Y_{n_s}\|^q+\|\bar{Y}_s\|^q\right\|_{L^{2p}(\Omega;\mathbb{R})} \cdot\left\|\bar{Y}_s-Y_{n_s}\right\|^2_{L^{4p}(\Omega;\mathbb{R}^d)}
\\ \leq & C_{p,T} h,
\end{split}
\end{equation}
where the polynomial growth condition (\ref{OLC6}) on $\mu''(x)$, Lemma \ref{yd}, H\"{o}lder's inequality and Jensen's inequality were also used.
Now we return to $\|\tilde{R}_1(\mu)\|_{L^{p}(\Omega;\mathbb{R}^d)}$. Replacing $f$ in \eqref{RT1} with $\mu$ gives
\begin{align}\label{R1p}
%\begin{split}
& \|\tilde{R}_1(\mu) \|_{L^p(\Omega;\mathbb{R}^d)}
\\ \leq &
\left\|\mu'(Y_{n_s}) (s-t_{n_s}) \tilde{\mu}(Y_{n_s}) \right\|_{L^p(\Omega;\mathbb{R}^d)}
+ \Big\|\mu'(Y_{n_s})\sum^m_{j_1,j_2=1}L^{j_1} \sigma_{j_2}(Y_{n_s})I_{j_1,j_2}^{t_{n_s},s}\Big\|_{L^p(\Omega;\mathbb{R}^d)} +\left\|R_1(\mu)\right\|_{L^p(\Omega;\mathbb{R}^d)}
\nonumber \\ \leq& h\left\|\mu'(Y_{n_s})\mu(Y_{n_s})\right\|_{L^{p}(\Omega;\mathbb{R}^d)}
+ \frac{1}{2}\sum^m_{j_1,j_2=1}\left\|\mu'(Y_{n_s})\Big(L^{j_1} \sigma_{j_2}(Y_{n_s})(\Delta W^{j_1}_s \Delta W^{j_2}_s-\delta_{j_1,j_2}h)\Big)\right\|_{L^p(\Omega;\mathbb{R}^d)}
\nonumber \\
 &+\left\|R_1(\mu)\right\|_{L^p(\Omega;\mathbb{R}^d)}.
\nonumber
%\end{split}
\end{align}
Following the same line as estimating (\ref{mb4}) and noticing that $s-t_{n_s} \leq h$, one can similarly arrive at
\begin{equation}\label{eq:DeltaW2}
\frac{1}{2}\sum_{j_1,j_2=1}^{m} \|\Delta W^{j_1}_s \Delta W^{j_2}_s-\delta_{j_1,j_2}h\|_{L^p(\Omega;\mathbb{R})} \leq c_{p,m}h.
\end{equation}
Further, using H\"{o}lder's inequality we derive from Lemma \ref{mblem0} that for $0\leq s \leq t \leq T, 1\leq j_1,j_2\leq m$
\begin{equation}\label{con8}
\begin{split}
 \|\mu'(Y_{n_s})\mu(Y_{n_s})\|_{L^{p}(\Omega;\mathbb{R}^d)}\leq \|\mu'(Y_{n_s})\|_{L^{2p}(\Omega;\mathbb{R}^{d\times d})}\cdot\|\mu(Y_{n_s})\|_{L^{2p}(\Omega;\mathbb{R}^d)}
< \infty,\\
\left\|\mu'(Y_{n_s})L^{j_1} \sigma_{j_2}(Y_{n_s})\right\|_{L^p(\Omega;\mathbb{R}^d)} \leq \|\mu'(Y_{n_s})\|_{L^{2p}(\Omega;\mathbb{R}^{d\times d})}\cdot \left\|L^{j_1} \sigma_{j_2}(Y_{n_s})\right\|_{L^{2p}(\Omega;\mathbb{R}^d)} < \infty.
\end{split}
\end{equation}
%$$
%\leq h\left\|\mu'(Y_{n_s})\mu(Y_{n_s})\right\|_{L^{p}(\Omega;\mathbb{R}^d)}
%+ \sup_{1\leq j_1,j_2\leq m}\left\|\mu'(Y_{n_s})L^{j_1} \sigma_{j_2}(Y_{n_s})\right\|_{L^p(\Omega;\mathbb{R}^d)}\cdot c_{p,m}h +\left\|R_1(\mu)\right\|_{L^p(\Omega;\mathbb{R}^d)}
%$$
Now, using the  independence of $Y_{n_s}$ and
$\Delta W_s^{j_1}, \Delta W_s^{j_2}$, and
combining (\ref{R1mu}), (\ref{R1p}), \eqref{eq:DeltaW2} and (\ref{con8}), one can show
\begin{equation}\label{R1pmu}
\|\tilde{R}_1(\mu)\|_{L^p(\Omega;\mathbb{R}^d)} \leq C_{p,T} h.
\end{equation}
In a similar way as estimating $\|\tilde{R}_1(\mu)\|_{L^p(\Omega;\mathbb{R}^d)}$, one can derive that
\begin{equation}\label{R1pf}
\|\tilde{R}_1(\sigma_{i})\|_{L^p(\Omega;\mathbb{R}^d)} \leq C_{p,T} h. \quad \square
\end{equation}
Our proof of Theorem \ref{conthm} also needs the following Burkholder-Davis-Gundy type inequality for discrete-time martingale (see Theorem 3.28 in \cite{KS91} and Lemma 5.1 in \cite{HJ09b}).
\begin{lem}\label{BDG}
Let $Z_1,...,Z_N: \Omega \rightarrow \mathbb{R}$ be $\mathcal {F}/\mathcal {B}(\mathbb{R})$-measurable mappings with $\mathbb{E}\|Z_n\|^p< \infty$ for all $n \in \{1,...,N\}$ and with $\mathbb{E} [Z_{n+1}|Z_1,...,Z_n] = 0$ for all $n \in \{1,...,N\}$. Then
\begin{equation}
\|Z_1 + ... + Z_n\|_{L^p(\Omega;\mathbb{R})} \leq c_{p}\left(\|Z_1\|^2_{L^p(\Omega;\mathbb{R})}+ ... +  \|Z_n\|^2_{L^p(\Omega;\mathbb{R})} \right)^{\frac{1}{2}}
\end{equation}
for every $p \in [2,\infty)$, where $c_p$ are constants dependent of $p$, but independent of $n$.
\end{lem}

{\it Proof of Theorem \ref{conthm}.}  Applying It\^{o}'s formula to \eqref{xyd} gives
\begin{equation}\label{eq:ito2}
\begin{split}
\|X_s-\bar{Y}_s\|^2 = & 2 \int_0^s \Big\langle X_u - \bar{Y}_u,  \mu(X_u)-\tilde{\mu}(Y_{n_u}) \Big\rangle du
\\ &
+ 2 \sum_{i=1}^m \int_0^s \Big\langle X_u - \bar{Y}_u,  \sigma_i(X_u) - \sigma_i(Y_{n_u})-\sum_{j=1}^m L^j \sigma_i(Y_{n_u}) \Delta W_u^j \Big\rangle d W_u^i
\\ &
+ \sum_{i=1}^m \int_0^s \Big\| \sigma_i(X_u) - \sigma_i(Y_{n_u})-\sum_{j=1}^m L^j \sigma_i(Y_{n_u}) \Delta W_u^j \Big\|^2 du.
\end{split}
\end{equation}
For the integrand of the first term in \eqref{eq:ito2}, we use \eqref{OLC1} and the Cauchy-Schwarz inequality to arrive at
\begin{eqnarray}
%\begin{split}
&&\Big\langle X_u - \bar{Y}_u,  \mu(X_u)-\tilde{\mu}(Y_{n_u}) \Big\rangle
\nonumber \\ &= &
\Big\langle X_u - \bar{Y}_u,  \mu(X_u)-\mu(\bar{Y}_{u}) \Big\rangle
+ \Big\langle X_u - \bar{Y}_u,  \mu(\bar{Y}_{u}) - \mu(Y_{n_u}) \Big\rangle
+  \Big\langle X_u - \bar{Y}_u, h\,\|\mu(Y_{n_u})\| \,\tilde{\mu}(Y_{n_u}) \Big\rangle
\nonumber \\
&\leq & K \|X_u - \bar{Y}_u\|^2 + \Big\langle X_u - \bar{Y}_u,  \mu(\bar{Y}_{u}) - \mu(Y_{n_u}) \Big\rangle
+ \frac{1}{2} \|X_u - \bar{Y}_u\|^2 + \frac{1}{2} h^2 \|\mu(Y_{n_u})\|^2\, \|\tilde{\mu}(Y_{n_u})\|^2
\nonumber \\
&\leq & (K+\frac{1}{2}) \|X_u - \bar{Y}_u\|^2 + \Big\langle X_u - \bar{Y}_u,
\mu(\bar{Y}_{u}) - \mu(Y_{n_u}) \Big\rangle + \frac{1}{2} h^2 \|\mu(Y_{n_u})\|^4.
%\end{split}
\label{eq:integrand1}
\end{eqnarray}
For the integrand of the third term in \eqref{eq:ito2}, one can use an elementary inequality, the notation \eqref{eq:sigmaTaylor} and \eqref{OLC3} to get
\begin{equation}\label{eq:integrand3}
\begin{split}
\Big\| \sigma_i(X_u) - \sigma_i(Y_{n_u})-\sum_{j=1}^m L^j \sigma_i(Y_{n_u}) \Delta W_u^j \Big\|^2 \leq&
2 \|\sigma_i(X_u) - \sigma_i(\bar{Y}_u)\|^2 + 2 \| \tilde{R_1}(\sigma_i) \|^2
\\ \leq&
2 K^2 \|X_u - \bar{Y}_u\|^2 + 2 \|\tilde{R_1}(\sigma_i)\|^2.
\end{split}
\end{equation}
Inserting \eqref{eq:integrand1} and \eqref{eq:integrand3} into \eqref{eq:ito2} and using the notation \eqref{eq:sigmaTaylor} show
\begin{equation}
\begin{split}
\|X_s-\bar{Y}_s\|^2 = & (2K + 1 + 2mK^2)\int_0^s \|X_u-\bar{Y}_u\|^2 du + h^2 \, \int_0^s \|\mu(Y_{n_u})\|^4 du
\\ & +
2 \sum_{i=1}^m \int_0^s \|\tilde{R_1}(\sigma_i)\|^2 du
+2\int_0^s \Big\langle X_u - \bar{Y}_u,  \mu(\bar{Y}_u)-\mu(Y_{n_u}) \Big\rangle du
\\ &
+ 2 \sum_{i=1}^m \int_0^s \Big\langle X_u - \bar{Y}_u,  \sigma_i(X_u) - \sigma_i(\bar{Y}_{u}) + \tilde{R_1}(\sigma_i) \Big\rangle d W_u^i.
\end{split}
\end{equation}
Hence
\begin{equation}\nonumber
\begin{split}
\sup_{0\leq s \leq t}\|X_s-\bar{Y}_s\|^2 \leq & (2K + 1 + 2mK^2)\int_0^t \|X_s-\bar{Y}_s\|^2 ds + h^2 \, \int_0^t \|\mu(Y_{n_s})\|^4 ds
\\ & +
2 \sum_{i=1}^m \int_0^t \|\tilde{R_1}(\sigma_i)\|^2 ds
+2 \sup_{0\leq s \leq t} \int_0^s \Big\langle X_u - \bar{Y}_u,  \mu(\bar{Y}_u)-\mu(Y_{n_u}) \Big\rangle du
\\ &
+ 2 \sup_{0\leq s \leq t} \sum_{i=1}^m \int_0^s \Big\langle X_u - \bar{Y}_u,  \sigma_i(X_u) - \sigma_i(\bar{Y}_{u}) + \tilde{R_1}(\sigma_i) \Big\rangle d W_u^i,
\end{split}
\end{equation}
and thus for $p \geq 4$
\begin{equation}\label{eq:expectation}
\begin{split}
& \Big\|\sup_{0\leq s \leq t}\|X_s-\bar{Y}_s\| \Big\|^2_{L^{p}(\Omega; \mathbb{R})} = \Big\|\sup_{0\leq s \leq t}\|X_s-\bar{Y}_s\|^2 \Big\|_{L^{\frac{p}{2}}(\Omega; \mathbb{R})}
\\ \leq &
(2K + 1 + 2mK^2)\int_0^t \|X_s-\bar{Y}_s\|_{L^{p}(\Omega; \mathbb{R}^d)}^2 ds + h^2 \, \int_0^t \|\mu(Y_{n_s})\|_{L^{2p}(\Omega; \mathbb{R}^d)}^4 ds
\\ & +
2 \sum_{i=1}^m \int_0^t \|\tilde{R_1}(\sigma_i)\|_{L^{p}(\Omega; \mathbb{R}^d)}^2 ds
+2 \Big\|\sup_{0\leq s \leq t} \int_0^s \Big\langle X_u - \bar{Y}_u,  \mu(\bar{Y}_u)-\mu(Y_{n_u}) \Big\rangle du \Big\|_{L^{\frac{p}{2}}(\Omega; \mathbb{R})}
\\ &
+ 2 \Big\|\sup_{0\leq s \leq t} \sum_{i=1}^m \int_0^s \Big\langle X_u - \bar{Y}_u,  \sigma_i(X_u) - \sigma_i(\bar{Y}_{u}) + \tilde{R_1}(\sigma_i) \Big\rangle dW_u^i\Big\|_{L^{\frac{p}{2}}(\Omega; \mathbb{R})}.
\end{split}
\end{equation}
For the last term on the right-hand side of \eqref{eq:expectation}, we use \eqref{OLC3}, \eqref{BDGc}, the Cauchy-Schwarz inequality and elementary inequalities to derive
\begin{eqnarray}
%\begin{split}
&&2 \Big\|\sup_{0\leq s \leq t} \sum_{i=1}^m \int_0^s \Big\langle X_u - \bar{Y}_u,  \sigma_i(X_u) - \sigma_i(\bar{Y}_{u}) + \tilde{R_1}(\sigma_i) \Big\rangle d W_u^i\Big\|_{L^{\frac{p}{2}}(\Omega; \mathbb{R})}
\nonumber \\ &\leq &
2\sum_{i=1}^m \Big\| \sup_{0\leq s \leq t} \int_0^s \Big\langle X_u - \bar{Y}_u,  \sigma_i(X_u) - \sigma_i(\bar{Y}_u) + \tilde{R}_1(\sigma_i) \Big\rangle d W_u^i\Big\|_{L^{\frac{p}{2}}(\Omega; \mathbb{R})}
\nonumber \\ &\leq &
p \sum_{i=1}^m \Big( \int_0^t \Big\|\Big\langle X_s - \bar{Y}_s,  \sigma_i(X_s) - \sigma_i(\bar{Y}_s) + \tilde{R}_1(\sigma_i) \Big\rangle \Big\|^2_{L^{\frac{p}{2}}(\Omega; \mathbb{R})} ds \Big)^{\frac{1}{2}}
\nonumber \\ &\leq &
p \sum_{i=1}^m \Big( \int_0^t \| X_s - \bar{Y}_s\|^2_{L^{p}(\Omega; \mathbb{R}^d)} \cdot \Big\| \sigma_i(X_s) - \sigma_i(\bar{Y}_s) + \tilde{R}_1(\sigma_i)\Big\|^2_{L^{p}(\Omega; \mathbb{R}^d)} ds \Big)^{\frac{1}{2}}
\nonumber \\ &\leq &
\sup_{0\leq s\leq t}\|X_s - \bar{Y}_s\|_{L^{p}(\Omega; \mathbb{R}^d)} \cdot p \sum_{i=1}^m \Big( \int_0^t \Big\| \sigma_i(X_s) - \sigma_i(\bar{Y}_s) + \tilde{R}_1(\sigma_i) \Big\|^2_{L^{p}(\Omega; \mathbb{R}^d)} ds  \Big)^{\frac{1}{2}}
\nonumber \\ &\leq &
\frac{1}{4}\sup_{0\leq s\leq t}\|X_s - \bar{Y}_s\|^2_{L^{p}(\Omega; \mathbb{R}^d)} + p^2 m \sum_{i=1}^m \int_0^t \Big\|\sigma_i(X_s) - \sigma_i(\bar{Y}_s) + \tilde{R}_1(\sigma_i) \Big\|^2_{L^{p}(\Omega; \mathbb{R}^d)} ds
\nonumber \\ &\leq &
\frac{1}{4}\sup_{0\leq s\leq t}\|X_s - \bar{Y}_s\|^2_{L^{p}(\Omega; \mathbb{R}^d)} + 2p^2 m^2 K \int_0^t \|X_s - \bar{Y}_s\|^2_{L^{p}(\Omega; \mathbb{R}^d)} ds
\nonumber\\ &&
+ 2p^2m
\sum_{i=1}^m \int_0^t \| \tilde{R}_1(\sigma_i) \|^2_{L^{p}(\Omega; \mathbb{R}^d)} ds.
%\end{split}
\label{eq:last}
\end{eqnarray}

At the same time, replacing $f$ in \eqref{taylorformula2} by $\mu$ and using the Cauchy-Schwarz inequality and elementary inequalities give
\begin{equation}\label{eq:lastbut1}
\begin{split}
&2 \Big\|\sup_{0\leq s \leq t} \int_0^s \Big\langle X_u - \bar{Y}_u,  \mu(\bar{Y}_u)-\mu(Y_{n_u}) \Big\rangle du \Big\|_{L^{\frac{p}{2}}(\Omega; \mathbb{R})}
\\ =&
2 \Big\|\sup_{0\leq s \leq t} \int_0^s \Big\langle X_u - \bar{Y}_u,  \mu'(Y_{n_u})\big(\sigma(Y_{n_u})\Delta W_u\big) + \tilde{R}_1(\mu) \Big\rangle du \Big\|_{L^{\frac{p}{2}}(\Omega; \mathbb{R})}
\\ \leq &
J +
2 \Big\|\sup_{0\leq s \leq t} \int_0^s \Big\langle X_u - \bar{Y}_u,  \tilde{R}_1(\mu) \Big\rangle du \Big\|_{L^{\frac{p}{2}}(\Omega; \mathbb{R})}
\\ \leq & J +
\int_0^t \|X_s-\bar{Y}_s\|^2_{L^{p}(\Omega; \mathbb{R}^d)} ds + \int_0^t \|\tilde{R}_1(\mu)\|^2_{L^{p}(\Omega; \mathbb{R}^d)} ds,
\end{split}
\end{equation}
where we denote
\begin{equation}\label{eq:J}
J = 2 \Big\|\sup_{0\leq s \leq t} \int_0^s \Big\langle X_u - \bar{Y}_u,  \mu'(Y_{n_u})\sigma(Y_{n_u})\Delta W_u \Big\rangle du \Big\|_{L^{\frac{p}{2}}(\Omega; \mathbb{R})}.
\end{equation}
Inserting \eqref{eq:last} and \eqref{eq:lastbut1} to \eqref{eq:expectation} yields
\begin{equation}\label{eq:Overall}
\begin{split}
&\frac{3}{4} \Big\|\sup_{0\leq s \leq t}\|X_s-\bar{Y}_s\| \Big\|^2_{L^{p}(\Omega; \mathbb{R})}
\\ \leq &
2(K + 1 + mK^2 + p^2 m^2K)\int_0^t \|X_s-\bar{Y}_s\|_{L^{p}(\Omega; \mathbb{R}^d)}^2 ds + h^2 \, \int_0^t \|\mu(Y_{n_s})\|_{L^{2p}(\Omega; \mathbb{R}^d)}^4 ds
\\ &
+ 2(1+p^2m)
\sum_{i=1}^m \int_0^t \| \tilde{R}_1(\sigma_i) \|^2_{L^{p}(\Omega; \mathbb{R}^d)} ds + \int_0^t \|\tilde{R}_1(\mu)\|^2_{L^{p}(\Omega; \mathbb{R}^d)} ds + J.
\end{split}
\end{equation}
Therefore it remains to estimate $J$ as \eqref{eq:J}. Using \eqref{Milstein2}, \eqref{eq:dx} and \eqref{eq:sigmaTaylor} shows
\begin{equation}\nonumber
\begin{split}
X_u - \bar{Y}_u =& X_{t_{n_u}} - Y_{n_u} + \int_{t_{n_u}}^u \mu(X_r) dr - \int_{t_{n_u}}^u \tilde{\mu}(Y_{n_u}) dr
\\ &+
\sum_{i=1}^m \int_{t_{n_u}}^u \Big[\sigma_{i}(X_r)-\sigma_{i}(Y_{n_r})-\sum_{j=1}^{m}L^{j}\sigma_{i}(Y_{n_r})\Delta W_{r}^{j} \Big] dW^i_r
\\ = &
\int_{t_{n_u}}^u \big[\mu(X_r)-\mu(X_{t_{n_u}})\big] dr + \sum_{i=1}^m \int_{t_{n_u}}^u \big[ \sigma_{i}(X_r)- \sigma_{i}(\bar{Y}_r) \big] dW^i_r + \sum_{i=1}^m \int_{t_{n_u}}^u \tilde{R}_1(\sigma_i)  dW^i_r
\\ &+
(u-t_{n_u}) \,\mu(X_{t_{n_u}}) - (u-t_{n_u}) \,\tilde{\mu}(Y_{n_u}) + X_{t_{n_u}} - Y_{n_u}.
\end{split}
\end{equation}
Thus
\begin{align}
%\begin{split}
J \leq & 2 \Big\|\sup_{0\leq s \leq t} \int_0^s \Big\langle \int_{t_{n_u}}^u \big[\mu(X_r)-\mu(X_{t_{n_u}})\big] dr,  \mu'(Y_{n_u})\sigma(Y_{n_u})\Delta W_u \Big\rangle du \Big\|_{L^{\frac{p}{2}}(\Omega; \mathbb{R})}
\nonumber \\ & +
2 \Big\|\sup_{0\leq s \leq t} \int_0^s \Big\langle \sum_{i=1}^m \int_{t_{n_u}}^u \big[ \sigma_{i}(X_r)- \sigma_{i}(\bar{Y}_r)\big]  dW^i_r,  \mu'(Y_{n_u})\sigma(Y_{n_u})\Delta W_u \Big\rangle du \Big\|_{L^{\frac{p}{2}}(\Omega; \mathbb{R})}
\nonumber \\ & +
2 \Big\|\sup_{0\leq s \leq t} \int_0^s \Big\langle \sum_{i=1}^m \int_{t_{n_u}}^u \tilde{R}_1(\sigma_i)  dW^i_r,  \mu'(Y_{n_u})\sigma(Y_{n_u})\Delta W_u \Big\rangle du \Big\|_{L^{\frac{p}{2}}(\Omega; \mathbb{R})}
\nonumber \\ & +
2 \Big\|\sup_{0\leq s \leq t} \int_0^s \Big\langle \zeta_{n_u},  \mu'(Y_{n_u})\sigma(Y_{n_u})\Delta W_u \Big\rangle du \Big\|_{L^{\frac{p}{2}}(\Omega; \mathbb{R})}
\nonumber \\ & +
2 \Big\|\sup_{0\leq s \leq t} \int_0^s \Big\langle X_{t_{n_u}} - Y_{n_u},  \mu'(Y_{n_u})\sigma(Y_{n_u})\Delta W_u \Big\rangle du \Big\|_{L^{\frac{p}{2}}(\Omega; \mathbb{R})}
\nonumber \\ := & J_{1} + J_{2} + J_{3} + J_{4} + J_{5},
%\end{split}
\label{eq:BJ}
\end{align}
where $\zeta_{n_u} \in \mathcal{F}_{t_{n_u}}$ is defined by
\begin{equation}\label{eq:zeta}
\zeta_{n_u} = (u-t_{n_u}) \,\mu(X_{t_{n_u}}) - (u-t_{n_u}) \,\tilde{\mu}(Y_{n_u}).
\end{equation}

To begin with, we establish the estimate
\begin{equation}\label{eq:frequent_ineq}
\|\mu'(Y_{n_u})\sigma(Y_{n_u})\Delta W_u \|_{L^{p}(\Omega; \mathbb{R}^d)} \leq C_{p,T} h^{\frac{1}{2}},
\end{equation}
which is frequently used later. Using H\"{o}lder's inequality, Lemma \ref{mblem0} and Lemma \ref{BDG-Continuous}, we know that for $0\leq k\leq N-1, t_k \leq s < t_{k+1}$
\begin{equation}\label{supcon04}
\left\| \mu'(Y_{k})\sigma(Y_{k})\right\|_{L^{2p}(\Omega;\mathbb{R}^{d\times m})} \leq \| \mu'(Y_{k})\|_{L^{4p}(\Omega;\mathbb{R}^{d\times d})}\|\sigma(Y_{k})\|_{L^{4p}(\Omega;\mathbb{R}^{d\times m})} < \infty
\end{equation}
and
\begin{equation}\label{supcon05}
\left\| W_s-W_{t_{k}} \right\|_{L^{2p}(\Omega;\mathbb{R}^m)} = \left\| W_{s-t_k}-W_{t_{0}} \right\|_{L^{2p}(\Omega;\mathbb{R}^m)} \leq 2p \sqrt{m} h^{1/2}.
\end{equation}
Combining \eqref{supcon04} and \eqref{supcon05} one can readily obtain \eqref{eq:frequent_ineq} by H\"{o}lder's inequality. Concerning $J_1$, we use H\"{o}lder's inequality, \eqref{d6} and \eqref{eq:frequent_ineq} to arrive at
\begin{equation}\label{eq:J1}
\begin{split}
J_1 \leq  2 \int_0^t \int_{t_{n_u}}^u \Big\|\mu(X_r) - \mu(X_{t_{n_u}})\Big\|_{L^{p}(\Omega; \mathbb{R}^d)} \cdot \Big\|\mu'(Y_{n_u})\sigma(Y_{n_u})\Delta W_u \Big\|_{L^{p}(\Omega; \mathbb{R}^d)} dr du \leq  C_{p,T} h^2.
\end{split}
\end{equation}
For $J_{2}$, using H\"{o}lder's inequality, \eqref{OLC3}, \eqref{BDGc}, elementary inequalities and \eqref{eq:frequent_ineq} gives
\begin{eqnarray}
%\begin{split}
J_{2} &\leq &2 \int_0^t \Big\|\sum_{i=1}^m \int_{t_{n_u}}^u \big[ \sigma_{i}(X_r)- \sigma_{i}(\bar{Y}_r)\big]  dW^i_r \Big\|_{L^{p}(\Omega; \mathbb{R}^d)} \cdot \Big\| \mu'(Y_{n_u})\sigma(Y_{n_u})\Delta W_u \Big\|_{L^{p}(\Omega; \mathbb{R}^d)} du
\nonumber \\ &\leq &
\int_0^t \frac{1}{h}\,\Big\|\sum_{i=1}^m \int_{t_{n_u}}^u \big[ \sigma_{i}(X_r)- \sigma_{i}(\bar{Y}_r) \big] dW^i_r \Big\|^2_{L^{p}(\Omega; \mathbb{R}^d)} du
+
\int_0^t h\, \Big\| \mu'(Y_{n_u})\sigma(Y_{n_u})\Delta W_u \Big\|^2_{L^{p}(\Omega; \mathbb{R}^d)} du
\nonumber \\ &\leq &
\frac{p^2}{h}\, \int_0^t \int_{t_{n_u}}^u \sum_{i=1}^m  \Big\| \sigma_{i}(X_r)- \sigma_{i}(\bar{Y}_r) \Big\|^2_{L^{p}(\Omega; \mathbb{R}^d)}dr du + C_{p,T} h^2
\nonumber \\ &\leq &
mp^2 K^2\, \int_0^t \sup_{0\leq r\leq u}  \big\| X_r - \bar{Y}_r \big\|^2_{L^{p}(\Omega; \mathbb{R}^d)} du + C_{p,T} h^2.
%\end{split}
\label{eq:J2}
\end{eqnarray}
For $J_{3}$, similarly as above we obtain that
\begin{equation}\label{eq:J3}
\begin{split}
J_{3} \leq & 2 \int_0^t \Big\|\sum_{i=1}^m \int_{t_{n_u}}^u \tilde{R}_1(\sigma_i)  dW^i_r \Big\|_{L^{p}(\Omega; \mathbb{R}^d)} \cdot \Big\| \mu'(Y_{n_u})\sigma(Y_{n_u})\Delta W_u \Big\|_{L^{p}(\Omega; \mathbb{R}^d)} du
\\ \leq & 2p \int_0^t \Big( \int_{t_{n_u}}^u \sum_{i=1}^m \Big\| \tilde{R}_1(\sigma_i)\Big\|^2_{L^{p}(\Omega; \mathbb{R}^d)} dr \Big)^{\frac{1}{2}} \cdot \Big\| \mu'(Y_{n_u})\sigma(Y_{n_u})\Delta W_u \Big\|_{L^{p}(\Omega; \mathbb{R}^d)} du
\\ \leq & C_{p,T} h^2,
\end{split}
\end{equation}
where \eqref{R_estimate} and \eqref{eq:frequent_ineq} were also used.
Now, it remains to estimate $J_{4}$ and $J_{5}$. We split $J_{4}$ into two terms as follows:
\begin{equation}\label{eq:J4}
\begin{split}
J_{4} \leq & 2 \left\|\sup_{0\leq s \leq t} \Big| \sum_{k=0}^{n_s-1} \int_{t_k}^{t_{k+1}} \Big\langle \zeta_{k},  \mu'(Y_{k})\sigma(Y_{k})\Delta W_u \Big\rangle du \Big| \right\|_{L^{\frac{p}{2}}(\Omega; \mathbb{R})}
\\ & + 2 \left\|\sup_{0\leq s \leq t} \Big|  \int_{t_{n_s}}^{s} \Big\langle \zeta_{n_s},  \mu'(Y_{n_s})\sigma(Y_{n_s})\Delta W_u \Big\rangle du \Big| \right\|_{L^{\frac{p}{2}}(\Omega; \mathbb{R})}
\\  := & J_{41} + J_{42}.
\end{split}
\end{equation}

Recall that $\zeta_{k} \in \mathcal{F}_{t_{k}}$ for $k=0,1,...,N-1$. It can be readily verified that the discrete time process
$$
\chi_n:=\left\{\sum_{k=0}^{n-1}\int_{t_k}^{t_{k+1}} \Big\langle \zeta_{k},  \mu'(Y_{k})\sigma(Y_{k})\Delta W_u \Big\rangle du\right\}, \:\:n \in \{0,1,...,N\}
$$
is an $\{\mathcal {F}_{t_n}: 0\leq n \leq N\}$-martingale. With the aid of Doob's maximal inequality, Lemma \ref{BDG} and H\"{o}lder's inequality we obtain that for $p \geq 4$
\begin{eqnarray}
%\begin{split}
J_{41} &\leq&  \frac{2p}{p-2} \Big\| \sum_{k=0}^{n_t-1} \int_{t_k}^{t_{k+1}} \Big\langle \zeta_{k},  \mu'(Y_{k})\sigma(Y_{k})\Delta W_u \Big\rangle du \Big\|_{L^{\frac{p}{2}}(\Omega; \mathbb{R})}
\nonumber \\ &\leq &
\frac{2pc_{p/2}}{p-2} \Big(\sum_{k=0}^{n_t-1} \Big\| \int_{t_k}^{t_{k+1}} \Big\langle \zeta_{k},  \mu'(Y_{k})\sigma(Y_{k})\Delta W_u \Big\rangle du  \Big\|^2_{L^{\frac{p}{2}}(\Omega;\mathbb{R})}\Big)^{1/2}
\nonumber \\ &\leq &
\frac{2pc_{p/2}}{p-2} \Big(\sum_{k=0}^{n_t-1} h \int_{t_k}^{t_{k+1}}  \Big\| \Big\langle \zeta_{k},  \mu'(Y_{k})\sigma(Y_{k})\Delta W_u \Big\rangle \Big\|^2_{L^{\frac{p}{2}}(\Omega;\mathbb{R})} du\Big)^{1/2}
\nonumber \\ &\leq&
\frac{2pc_{p/2}}{p-2} \Big(\sum_{k=0}^{n_t-1} h \int_{t_k}^{t_{k+1}}
\big\| \zeta_{k}\big\|^2_{L^{p}(\Omega;\mathbb{R}^d)} \cdot \big\|\mu'(Y_{k})
\sigma(Y_{k})\Delta W_u \big\|^2_{L^{p}(\Omega;\mathbb{R}^d)} du\Big)^{1/2},
%\end{split}
\label{supcon4}
\end{eqnarray}
where by \eqref{eq:zeta}, Lemma \ref{mblem0} and Lemma \ref{mblem}
\begin{equation}\label{eq:zeta_estimate}
\big\| \zeta_{k}\big\|_{L^{p}(\Omega;\mathbb{R}^d)} \leq h \|\mu(X_{t_k})\|_{L^{p}(\Omega;\mathbb{R}^d)} + h \|\tilde{\mu}(Y_k)\|_{L^{p}(\Omega;\mathbb{R}^d)} \leq C_{p,T}h.
\end{equation}

Hence, taking \eqref{eq:frequent_ineq} and \eqref{eq:zeta_estimate} into account, we derive from (\ref{supcon4})  that
\begin{equation}\label{supcon005}
J_{41} \leq C_{p,T} h^2.
\end{equation}
For the second term $J_{42}$, H\"{o}lder's inequality, \eqref{eq:zeta_estimate} and \eqref{eq:frequent_ineq} give that for $p \geq 4$
\begin{equation}\label{supcon5}
\begin{split}
\Big(\frac{J_{42}}{2}\Big)^{\frac{p}{2}}=& \mathbb{E}\Big(\sup_{0\leq s\leq t}\Big|\int_{t_{n_s}}^{s} \langle \zeta_{n_s},  \mu'(Y_{n_s})\sigma(Y_{n_s})\Delta W_u \rangle du \Big|\Big)^{\frac{p}{2}}
\\ \leq &
h^{\frac{p}{2}-1} \mathbb{E} \Big( \sup_{0\leq s\leq t} \int_{t_{n_s}}^{s}\Big|\left\langle \zeta_{n_s},  \mu'(Y_{n_s})\sigma(Y_{n_s})\Delta W_u \right\rangle\Big|^{\frac{p}{2}} du \Big)
\\ \leq &
h^{\frac{p}{2}-1} \mathbb{E} \Big( \sum_{k=0}^{n_t-1} \int_{t_k}^{t_{k+1}}\Big|\langle \zeta_{k},  \mu'(Y_{k})\sigma(Y_{k})\Delta W_u \rangle\Big|^{\frac{p}{2}}du + \int_{t_{n_t}}^{t}\Big|\langle \zeta_{n_t},  \mu'(Y_{n_t})\sigma(Y_{n_t})\Delta W_u \rangle\Big|^{\frac{p}{2}}du \Big)
\\ = &
h^{\frac{p}{2}-1}  \int_{0}^{t} \mathbb{E} \Big|\langle \zeta_{n_u},  \mu'(Y_{n_u})\sigma(Y_{n_u})\Delta W_u \rangle\Big|^{\frac{p}{2}}du
\\ \leq &
h^{\frac{p}{2}-1}  \int_{0}^{t} \|\zeta_{n_u}\|^{\frac{p}{2}}_{L^p(\Omega;\mathbb{R}^d)} \cdot
\|\mu'(Y_{n_u})\sigma(Y_{n_u})\Delta W_u\|^{\frac{p}{2}}_{L^p(\Omega;\mathbb{R}^d)} du
\\ \leq & C_{p,T} h^{\frac{5p}{4}-1},
\end{split}
\end{equation}
which implies that for $p \geq 4$ and $h \in (0,1]$ there exists a suitable constant $C_{p,T}$ such that
\begin{equation} \label{supcon7}
J_{42} \leq 2C^{\frac{2}{p}}_{p,T} h^{\frac{5p-4}{2p}}  \leq C_{p,T} h^2.
\end{equation}
Gathering (\ref{supcon005}) and (\ref{supcon7}) we derive from \eqref{eq:J4} that for $p \geq 4$
\begin{equation}\label{supcon006}
J_4 \leq C_{p,T} h^2.
\end{equation}
Similarly, we split $J_5$ as follows:
\begin{equation}\label{eq:BJ5}
\begin{split}
J_{5} \leq & 2 \left\|\sup_{0\leq s \leq t} \Big| \sum_{k=0}^{n_s-1} \int_{t_k}^{t_{k+1}} \Big\langle X_{t_{k}} - Y_{k},  \mu'(Y_{k})\sigma(Y_{k})\Delta W_u \Big\rangle du \Big| \right\|_{L^{\frac{p}{2}}(\Omega; \mathbb{R})}
\\ & + 2 \left\|\sup_{0\leq s \leq t} \Big|  \int_{t_{n_s}}^{s} \Big\langle X_{t_{n_s}} - Y_{n_s},  \mu'(Y_{n_s})\sigma(Y_{n_s})\Delta W_u \Big\rangle du \Big| \right\|_{L^{\frac{p}{2}}(\Omega; \mathbb{R})}
\\  := & J_{51} + J_{52}.
\end{split}
\end{equation}
With regard to $J_{51}$, following the same line as \eqref{supcon4} yields
\begin{equation}\label{eq:J51}
\begin{split}
J_{51} \leq &
\frac{2pc_{p/2}}{p-2} \Big(\sum_{k=0}^{n_t-1} h \int_{t_k}^{t_{k+1}}  \| X_{t_{k}} - Y_{k}\|^2_{L^{p}(\Omega;\mathbb{R}^d)} \cdot \|\mu'(Y_{k})\sigma(Y_{k})\Delta W_u \|^2_{L^{p}(\Omega;\mathbb{R}^d)} du\Big)^{1/2}
\\ \leq &
\sup_{0\leq s \leq t} \| X_{s} - \bar{Y}_{s}\|_{L^{p}(\Omega;\mathbb{R}^d)} \cdot \frac{2pc_{p/2}}{p-2} \Big(\sum_{k=0}^{n_t-1} h \int_{t_k}^{t_{k+1}} \|\mu'(Y_{k})\sigma(Y_{k})\Delta W_u \|^2_{L^{p}(\Omega;\mathbb{R}^d)} du\Big)^{1/2}
\\ \leq &
\frac{1}{4} \sup_{0\leq s \leq t} \| X_{s} - \bar{Y}_{s}\|^2_{L^{p}(\Omega;\mathbb{R}^d)} + \frac{4p^2c_{p/2}^2}{(p-2)^2} \sum_{k=0}^{n_t-1} h \int_{t_k}^{t_{k+1}} \|\mu'(Y_{k})\sigma(Y_{k})\Delta W_u \|^2_{L^{p}(\Omega;\mathbb{R}^d)} du
\\ \leq &
\frac{1}{4} \sup_{0\leq s \leq t} \| X_{s} - \bar{Y}_{s}\|^2_{L^{p}(\Omega;\mathbb{R}^d)} + C_{p,T} h^2,
\end{split}
\end{equation}
where the fact $\bar{Y}_{t_k} = Y_k, k =0,1,...,N-1$,  elementary inequalities and \eqref{eq:frequent_ineq} were used.
For $J_{52}$, we follow the same way as \eqref{supcon5} to obtain
\begin{equation}\label{eq:BJ52}
\begin{split}
\Big( \frac{J_{52}}{2}\Big)^{\frac{p}{2}}\leq &
h^{\frac{p}{2}-1}  \int_{0}^{t} \|X_{t_{n_u}}-Y_{n_u}\|^{\frac{p}{2}}_{L^p(\Omega;\mathbb{R}^d)} \cdot
\|\mu'(Y_{n_u})\sigma(Y_{n_u})\Delta W_u\|^{\frac{p}{2}}_{L^p(\Omega;\mathbb{R}^d)} du
\\ \leq &
\sup_{0\leq s \leq t} \| X_{s} - \bar{Y}_{s}\|^{\frac{p}{2}}_{L^{p}(\Omega;\mathbb{R}^d)} \cdot C_{p,T} h^{\frac{3p}{4}-1}.
\end{split}
\end{equation}
Thus
\begin{equation}\label{eq:J52}
\begin{split}
J_{52} \leq 2 \sup_{0\leq s \leq t} \| X_{s} - \bar{Y}_{s}\|_{L^{p}(\Omega;\mathbb{R}^d)} \cdot C_{p,T}^{\frac{2}{p}} \, h^{\frac{3p-4}{2p}} \leq \frac{1}{4} \sup_{0\leq s \leq t} \| X_{s} - \bar{Y}_{s}\|^2_{L^{p}(\Omega;\mathbb{R}^d)} + 4 C_{p,T}^{\frac{4}{p}} \, h^{\frac{3p-4}{p}}.
\end{split}
\end{equation}
Note that $\frac{3p-4}{4} \geq 2$ for $p \geq 4$ and that $h \leq 1$. Plugging \eqref{eq:J51} and \eqref{eq:J52} into \eqref{eq:BJ5} yields
\begin{equation}\label{eq:J5}
J_5 \leq \frac{1}{2} \sup_{0\leq s \leq t} \| X_{s} - \bar{Y}_{s}\|^2_{L^{p}(\Omega;\mathbb{R}^d)} + C_{p,T} h^2.
\end{equation}
Inserting \eqref{eq:J1}, \eqref{eq:J2}, \eqref{eq:J3}, \eqref{supcon006} and \eqref{eq:J5} into \eqref{eq:BJ} leads to
\begin{equation}\label{eq:JF}
J \leq mp^2 K^2\, \int_0^t \sup_{0\leq r\leq u}  \Big\| X_r - \bar{Y}_r \Big\|^2_{L^{p}(\Omega; \mathbb{R}^d)} du + \frac{1}{2} \sup_{0\leq s\leq t}\Big\| X_s - \bar{Y}_s \Big\|^2_{L^{p}(\Omega; \mathbb{R}^d)}  + C_{p,T} h^2.
\end{equation}
Hence, using estimates in Lemma \ref{mblem0} and Lemma \ref{Rlem} we derive from \eqref{eq:Overall} that
\begin{equation}\nonumber
\begin{split}
\frac{3}{4} \Big\|\sup_{0\leq s \leq t}\|X_s-\bar{Y}_s\| \Big\|^2_{L^{p}(\Omega; \mathbb{R})}
\leq &
2(K + 1 + mK^2 + p^2 m^2K + \frac{1}{2}mp^2K^2)\int_0^t \sup_{0\leq u\leq s}\|X_u-\bar{Y}_u\|_{L^{p}(\Omega; \mathbb{R}^d)}^2 ds
\\ &
+ \frac{1}{2} \sup_{0\leq s\leq t}\Big\| X_s - \bar{Y}_s \Big\|^2_{L^{p}(\Omega; \mathbb{R}^d)}  + C_{p,T} h^2.
\end{split}
\end{equation}
Thus $\big\|\sup_{0\leq s \leq t}\|X_s-\bar{Y}_s\| \big\|^2_{L^{p}(\Omega; \mathbb{R})}$ is finite by Lemma \ref{mblem} and
\begin{equation}\label{eq:final}
\begin{split}
&\Big\|\sup_{0\leq s \leq t}\|X_s-\bar{Y}_s\| \Big\|^2_{L^{p}(\Omega; \mathbb{R})}
\leq
C_{p,T} \int_0^t \Big\|\sup_{0\leq u \leq s}\|X_u-\bar{Y}_u\| \Big\|^2_{L^{p}(\Omega; \mathbb{R})} ds + C_{p,T} h^2.
\end{split}
\end{equation}
The Gronwall inequality gives the desired result for $p \geq 4$. Using H\"{o}lder's inequality gives the assertion for $1\leq p <4$ and the proof is complete. $\square$

\section{An illustrative example}

\begin{figure}[h]
         \centering
         \includegraphics[width=3.5in,height=2.5in]{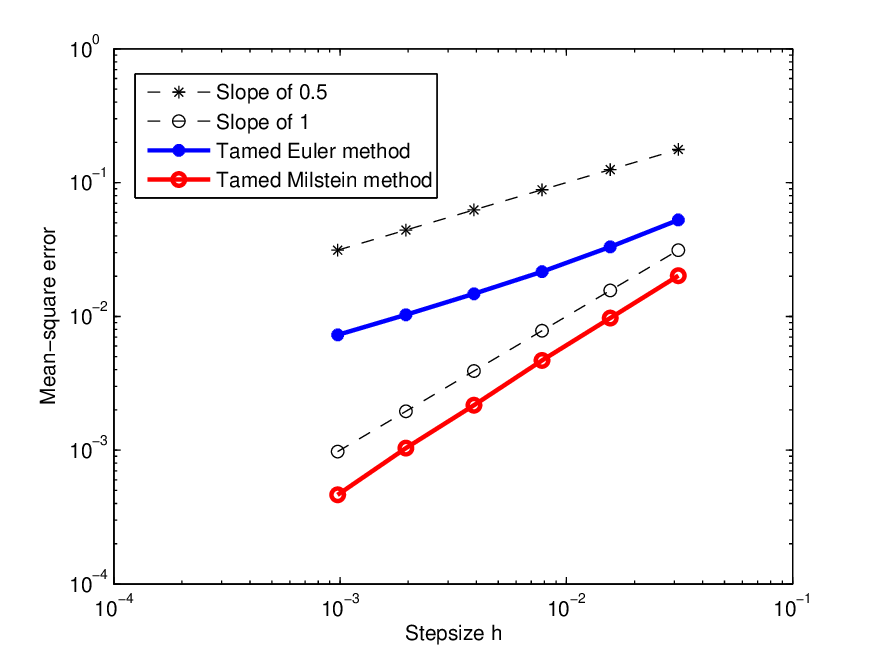}
         \caption{Root mean square approximation error versus  stepsize $h$
         to approximate (\ref{test}). }
         \label{fig1}
\end{figure}

\begin{figure}[h]
         \centering
         \includegraphics[width=3.5in,height=2.5in]{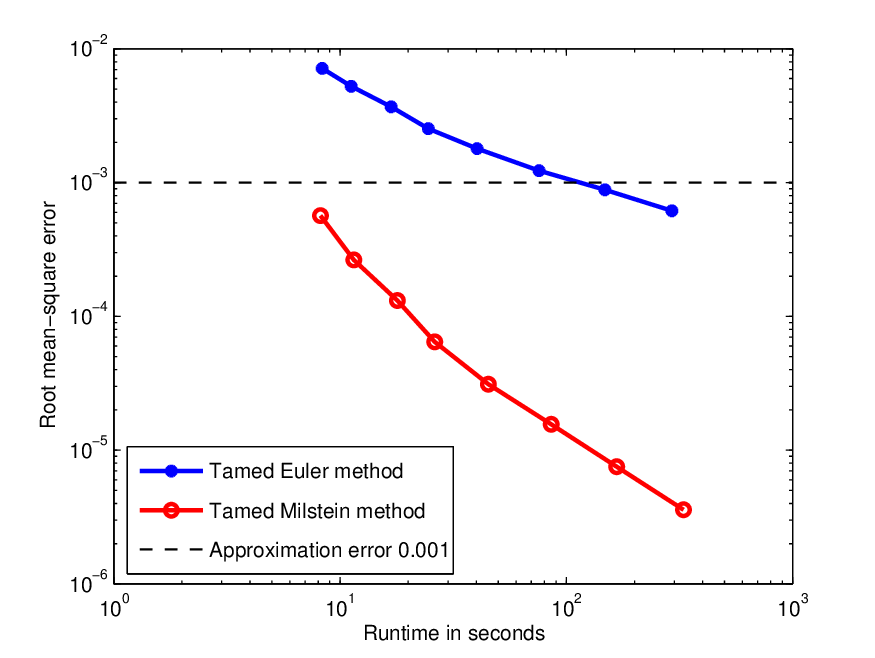}
         \caption{Root mean square approximation error versus runtime for $N \in \{2^{10},...,2^{17}\}$. }
         \label{fig2}
\end{figure}

In \cite{HJK10}, the authors have demonstrated the computational efficiency of the tamed Euler scheme, compared to the implicit Euler method. In this section we compare computational efficiency of the tamed Milstein scheme and the tamed Euler scheme. To this end we choose a simple SDE (\ref{sdes})
\begin{equation}\label{test}
d X_t = -X_t^5 dt + X_t d W_t, \quad X_0=1
\end{equation}
for $t \in [0,1]$.
Figure \ref{fig1} depicts the root mean-square errors (\ref{Err}) as a function of the stepsize $h$ in log-log plot, where the expectation is approximated by the mean of 5000 independent realizations. As expected, the tamed Milstein scheme gives an error that decreases proportional to $h$, whereas the tamed Euler scheme gives errors that decrease proportional to $h^{\frac{1}{2}}$.  To show the efficiency of the tamed Milstein method clearly, we present in Figure \ref{fig2} the root mean-square errors of both methods as function of the runtime when $N \in \{2^{10},...,2^{17}\}$ and the mean of 1000 independent  paths are used to approximate the expectation in (\ref{Err}). Suppose that the strong approximation problem (\ref{Err}) of the SDE (\ref{test}) should be solved with the precision $\varepsilon = 0.001$. From Figure \ref{fig2}, one can detect that $N=2^{10}$ in the case of the tamed Milstein method (\ref{MilsteinC}) and that $N=2^{16}$ in the case of the tamed Euler method (\ref{tamedEuler}) achieves the desired precision $\varepsilon = 0.001$ in (\ref{Err}). Moreover, the tamed Milstein scheme requires 8.1860 seconds while the tamed Euler scheme requires 147.9230 seconds to achieve the precision $\varepsilon = 0.001$ in (\ref{Err}). The tamed Milstein method is for the SDE (\ref{test}) with commutative noise thus much faster than the tamed Euler method.

\end{document}